\documentclass{article}
\usepackage{amsmath,amssymb}
\usepackage{graphicx}
\usepackage{float,afterpage}
\usepackage[active]{srcltx}
\usepackage{color}
\newcommand{\arcosh}{\mathop{\mathrm{arcosh}}}
\newtheorem{theorem}{Theorem}[section]
\newtheorem{rem}[theorem]{Remark}

\newtheorem{proposition}[theorem]{Proposition}
\newtheorem{cor}[theorem]{Corollary}
\newtheorem{exmp}[theorem]{Example}

\newenvironment{remark}{\begin{rem}\rm}{\end{rem}}

\newcommand{\diag}{\mathop{\mathrm{diag}}}
\renewcommand{\epsilon}{\varepsilon}
\renewcommand{\tilde}{\widetilde}
\renewcommand{\hat}{\widehat}

     \newcommand{\NN}{\mathbb{N}}
     \newcommand{\ZZ}{\mathbb{Z}}
     \newcommand{\bA}{\mathbf{A}}
     \newcommand{\bB}{\mathbf{B}}
     \newcommand{\bH}{\mathbf{H}}
     \newcommand{\bK}{\mathbf{K}}
     \newcommand{\bT}{\mathbf{T}}
     \newcommand{\bU}{\mathbf{U}}
     \newcommand{\bW}{\mathbf{W}}
     
\newcommand{\diverg}{\mathop{\mathrm{div}}}
\newcommand{\grad}{\mathop{\mathrm{grad}}}
\DeclareMathOperator{\dist}{dist}

\newcommand{\hm}[1]{\leavevmode{\marginpar{\tiny%
$\hbox to 0mm{\hspace*{-0.5mm}$\leftarrow$\hss}%
\vcenter{\vrule depth 0.1mm height 0.1mm width \the\marginparwidth}%
\hbox to
0mm{\hss$\rightarrow$\hspace*{-0.5mm}}$\\\relax\raggedright #1}}}

\title{Spectral gap estimates for some block matrices
}
\author{Ivan Veseli\'c\thanks{TU-Chemnitz,
Fakult\"{a}t f\"{u}r Mathematik, 09107 Chemnitz, Germany, \newline e-mail:
ivan.veselic@mathematik.tu-chemnitz.de}, Kre\v simir Veseli\'c
\thanks{Fernuniversit\" at Hagen,
Fakult\" at f\" ur Mathematik und Informatik,
Postfach 940, D-58084 Hagen, Germany, e-mail:
kresimir.veselic@fernuni-hagen.de}
}
\date{ }
\begin{document}
\maketitle
\begin{abstract}
We estimate the size of the spectral gap at zero for some 
Hermitian block 
matrices. Included are quasi-definite matrices, quasi-semidefinite matrices
 (the closure 
of the set of the
quasi-definite matrices) and some related block matrices which
need not belong to either of these classes.
Matrices of such structure arise in quantum models of 
possibly disordered systems
with supersymmetry or graphene like symmetry.
Some of the results immediately extend to infinite dimension.
\end{abstract}
\setcounter{equation}{0}
\setcounter{figure}{0}
\section{Introduction}\label{Introduction}

Consider (finite) Hermitian block matrices
of the form
\begin{equation}\label{quasidef}
H =
\left[\begin{array}{rr}
A     &       B \\
B^*   &      -C  \\
\end{array}\right]
\end{equation}
(the minus sign is set by convenience).
If \(A,C\)  are positive definite then the matrix \(H\) is called
quasi-definite. 
These matrices have several
remarkable properties, one of them being that they are
always nonsingular with a spectral gap at zero,
\begin{equation}\label{gapHH0}
\rho(H) \supseteq (-\min\sigma(C),\min\sigma(A))
\end{equation}
(\(\rho\) the resolvent set, \(\sigma\) the spectrum).
That is, the spectral gap of the block-diagonal part
of \(H\) in (\ref{quasidef}) can only grow if any
\(B\) is added. 
Moreover, quasi-definite matrices have two remarkable
{\em monotonicity properties}:
\begin{itemize}
\item[(A)] 
If \(B\) is replaced by \(tB\), \(t > 0\),
then all the eigenvalues go monotonically asunder
as \(t\) is growing (\cite{thompson},\cite{vankemp}).
\item[(B)] The same holds if \(A,C\) is replaced by
\(A + tI,C + tI\), \(t > 0\), respectively
(\cite{ikr_george}).
\end{itemize}
 
  In this note we study some related classes of matrices.
If in (\ref{quasidef})  the blocks \(A,C\) 
 are allowed to be only positive semidefinite
then \(H\)
will naturally be called {\em quasi-semidefinite}.
These matrices need not to be invertible.

It is relatively easy to characterise the nonsingularity of a quasi-semidefinite matrix, see Proposition
\ref{singularity} below. Giving
estimates for the gap at zero is more involved
and this note offers some results in this direction. 
Since the size of the spectral gap at zero is bounded from below by 
\(2\|H^{-1}\|^{-1}\) we will
give various bounds for this quantity in terms of the blocks
\(A,B,C\) where \(A,C\) are only positive semidefinite
and the properties of \(B\) come into play. It is known
from \cite{ves} that
the invertibility of \(B\) carries over to \(H\), but no
bound for \(H^{-1}\) was provided there. Some bounds for \(H^{-1}\)
were given in \cite{winklmeier}.

As a technical tool we derive a bound for the matrix \((I + AC)^{-1}\)
with \(A,C\) 
positive semidefinite
which might be of independent interest. We also sketch a 
related functional calculus for such products. 
More specifically, the present article provides the following.

\begin{enumerate}
\item A bound for \((I + AC)^{-1}\)
      with \(A,C\) positive semidefinite.
\item A characterisation of the nonsingularity of a quasi-semidefinite matrix.
\item A bound for \(H^{-1}\) based on the bound for \(B^{-1}\) including
an immediate generalisation to unbounded Hilbert-space operators defined
by quadratic forms. To this general environment we also extend an
elegant estimate obtained by \cite{kirsch} for the special case
\(A = C,\ B = B^*\).
\item A bound for \(H^{-1}\) based on the geometry of the null-spaces
of all of \(A,B,C\) and certain restrictions of these operators
to the orthogonal complements of these null-spaces.
\item Several counterexamples; some of them 
showing that some plausibly looking generalisations of the properties (A),
(B) above are not valid. 
\item A monotonicity and a sharp spectral inclusion result
for the case of Stokes matrices (those with \(C = 0\)).
\item A study of the spectral gap of a 
particular class of matrices which arise in the quantum mechanical
modelling of disordered systems (see e.g. \cite{Stolz}).
There we have \(C = A\) in (\ref{quasidef}), but \(A\) is not 
necessarily positive definite. In particular, we will illustrate
how changing boundary conditions can remove spurious eigenvalues
from the gap.
This is a specific, thoroughly worked out example on how to deal 
successfully with what is called spectral pollution.
\end{enumerate}

Let us remark that the variety of
special cases as well as techniques which we use illustrate the fact
that we did not succeed in obtaining a unified framework 
for spectral gap estimates for general quasi-definite matrices.

Quasi-semidefinite matrices and their infinite
dimensional analogs have important applications in Mathematical
Physics.
Although we here have no space to discuss the relevant 
models  in detail, we would like to convey an impression of the questions 
arising in this context. These have been the motivation for 
much of the research presented here.
Certain types of Dirac operators  are important examples.
In these cases the nonsingularity of \(H\) is typically due to the one
of \(B\) (see \cite{winklmeier} where this phenomenon was dubbed
'off-diagonal dominance').


Another particular motivation are quantum mechanical models of 
disordered solids. 
While this is a well established research field, recently
there has been interest in such models which give rise to
operators with block-structure, see e.g.~\cite{kirsch} or \cite{Stolz}.
For some of these models the block structure is a consequence of the Dirac-like symmetry arising in Hamiltonians describing graphene.     

Let us describe some of the specific spectral features which are of interest 
in this context.
We consider several instances of 
\emph{one-parameter Hermitian pencils}
 $A+t\, B, t\in \mathbb{R}$.
The well known \emph{monotonicity property}, namely that the eigenvalues
of  $A+t\, B$ grow monotonically in \(t\), if \(B\) is positive semidefinite
can, at least partly, be carried over to quasi-definite matrices as show
the properties (A), (B) listed above.
Here a question of particular importance is whether and how fast the spectral 
gap increases as $t$ grows. 
Several theorems of this paper provide answers to this question in 
specific situations.

As mentioned, certain physical models of 
disordered systems give rise to block-structured operator families.
In this context, estimates have to take into account the following
two important aspects.\\

(I)\ The size of the original physical system is macroscopic, i.e.~essentially infinite.
A mathematical understanding of the physical situation 
is -- as a rule  -- only possible 
by analysing larger and larger finite sample systems which describe the original physical situation
in the thermodynamic limit.

This leads to finite matrices or to operators with compact resolvent.
In any case, effectively one can reduce the focus on a finite number, say $n$, 
of eigenvalues, when analysing  monotonicity properties.
However, $n$ is not fixed but growing unboundedly as one passes to 
larger and larger sample scales.

Thus, efficient estimates on spectral gaps (or derivatives of eigenvalues)
are not allowed to depend on the system size -- expressed in 
the dimension of the matrix or the number of eigenvalues $n$.
We will pay special attention to this issue in the following.

(II)\ 
Due to the fact that one wants to model a disordered system, with a large number of degrees of freedom,
there is in fact not just \emph{one} coupling constant $ \in \mathbb{R}$, 
but rather a whole collection $(t_j)_{j\in \ZZ}$ of them. Thus the 
considered operator pencil is originally of the form 
\[
A+ \sum_j t_j \, B_j\ .
\]
A one-parameter family arises if one freezes all coupling constants except for one.
As a consequence, one is not dealing with one fixed unperturbed operator $A$, 
but rather with a whole collection of them, depending on the background 
configuration of the (other) coupling constants $(t_j, j \neq 0)$.
For this reason it would be desirable to obtain estimates on the spectral gap 
which do not depend on specific features of $A$.\\

The plan of the paper is as follows.
In the next section we provide certain basic preliminary estimates for quasi-semidefinite 
matrices. In Section \ref{quasi} the main results 
concerning the spectral gap size of such matrices are stated. 
These results are formulated for finite matrices.
In Section \ref{unbounded} we explain which results carry 
immediately over to the setting 
of (possibly unbounded) operators defined as quadratic forms. 
This includes the mentioned generalisation of an estimate from \cite{kirsch}, 
as well as a comparison with bounds obtained in \cite{winklmeier}.
In Section 
\ref{Stokes matrix}
we consider Stokes matrices. By reduction to a quadratic
eigenvalue problem we (i) prove monotonicity properties of the 
eigenvalues (but {\em not} as it would be naively expected from
cases (A) and (B) above), then (ii) give a tight bound for 
the two eigenvalues closest to zero.
The last section considers a special class of finite difference matrices, 
not necessarily quasi-definite, studied in \cite{Stolz}. 
Here we show that a stable spectral gap at zero can be achieved by 
an appropriate tuning of boundary conditions. Similar phenomena,
yet without rigorous proofs, are numerically observed on related models
with random diagonal entries.

\section{Some preliminary results} \label{preliminary}
To set the stage we collect some rather elementary statements and estimates.
\begin{proposition}\label{singularity}
A quasi-semidefinite matrix 
\[
H =
\left[\begin{array}{rr}
A     &       B \\
B^*   &      -C  \\
\end{array}\right].
\]
 is singular if and only if at least one of the subspaces
\[
\mathcal{N}(A)\cap \mathcal{N}(B^*),\quad
\mathcal{N}(C)\cap \mathcal{N}(B)
\]
{\rm(}\(\mathcal{N}\) denoting the null-space{\rm)}
is non-trivial. Moreover, in the obvious notation,
\begin{equation}\label{N(H)}
\mathcal{N}(H) =
\left[\begin{array}{cc}
\mathcal{N}(A)\cap \mathcal{N}(B^*) \\
\mathcal{N}(C)\cap \mathcal{N}(B)   \\
\end{array}\right].
\end{equation}
The value \(\min(\sigma(A))\) is an eigenvalue of \(H\) if and
only if \(\mathcal{N}(B^*)\) = \{0\} {\rm(}and similarly for \(\min(-\sigma(C))\){\rm)}.
\end{proposition}
{\bf Proof.} The equations
\[
Ax + By = 0,\quad B^*x - Cy = 0
\]
imply
\[
x^*Ax + x^*By = 0,\quad y^*B^*x - y^*Cy = 0.
\]
Since both \(x^*Ax\) and \(y^*Cy\) are real and non-negative, the
same is true of \(\pm x^*By\) such that, in fact, all three expressions
vanish. Since \(A,C\) are Hermitian positive semidefinite this implies
\(Ax = 0\) and \(Cy = 0\), then also \(B^*x = 0\) and \(By = 0\).
This proves (\ref{N(H)}); for the last assertion apply (\ref{N(H)})
to the matrix \(H - \min(\sigma(A)) I\). The other assertions follow trivially.
Q.E.D.\\

\begin{cor}\label{cor_singularity}
The matrix \(H\) is nonsingular if and only if the matrices
\[
A + \sqrt{BB^*},\quad C + \sqrt{B^*B}
\]
are positive definite.
\end{cor}
\begin{cor}\label{cor_parameter}
The null-space of the matrix \(H\) from (\ref{quasidef})
does not change if \(A\) is replaced by \(tA\), \(t > 0\)
{\rm }and similarly with \(B,C\){\rm }.
\end{cor}

To quantify the influence of \(B\) in (\ref{quasidef}) on 
the spectral gap in the quasi-definite case
we proceed as follows. First note the fundamental equality,
valid for all selfadjoint operators, saying that
\begin{equation}\label{equality}
\|(H - \lambda I)^{-1}\| 
=
\dist(\lambda,\sigma(H)).
\end{equation}
Taking any \(\lambda \) from the open interval
 \((-\min\sigma(C),\min\sigma(A))\) we have
\[
H - \lambda I = 
\left[
\begin{array}{ll}
(A - \lambda I)^{1/2} & 0  \\    
0 & (C + \lambda I)^{1/2}  \\
\end{array} 
\right]
W
\left[
\begin{array}{ll}
(A - \lambda I)^{1/2} & 0  \\    
0 & (C + \lambda I)^{1/2}  \\
\end{array} 
\right]
\]
with
\[
W = \left[
\begin{array}{ll}
I   &  Z  \\    
Z^* & -I  \\
\end{array} 
\right],
\quad
Z = (A - \lambda I)^{-1/2}B(C + \lambda I)^{-1/2}.
\]
As it is readily seen, the eigenvalues of the matrix \(W\)
are \(\pm \sqrt{1 + \sigma_i^2}\), where \( \sigma_i\)
are the singular values of \(Z\)
(cf.~\cite{saunders}).
Thus,
\[
\|W\| = \sqrt{\|I + Z^*Z\|},\quad 
\|W^{-1}\| = \sqrt{\|(I + Z^*Z)^{-1}\|}                .
\]
This gives the estimate
\begin{equation}\label{Hlambda}
\|(H - \lambda I)^{-1}\| 
\leq 
\frac{\sqrt{\|(I + Z^*Z)^{-1}\|}}
{\min\{\min(\sigma(A) - \lambda),\min(\sigma(C) + \lambda) }.
\end{equation}
Therefore by taking e.g.~\(\lambda = \lambda_0 =
\frac{1}{2}(\min(\sigma(A) - \min(\sigma(C))\)
we obtain

\begin{equation}\label{Hstretch}
\|(H - \lambda_0 I)^{-1}\| \leq 
\frac{2\sqrt{\|(I + Z^*Z)^{-1}}\|}{\min(\sigma(A) + \min(\sigma(C) }.
\end{equation}
We see that the gap is stretched at least by the factor
\[
\|(I + Z^*Z)^{-1}\|^{-1/2} = \sqrt{1 + \min_i\sigma_i^2}.
\]
\section{Spectral bounds for quasi-semidefinite matrices.}\label{quasi}
We will particularly be interested in how the appearance of the block
\(B\) can create a spectral gap at zero if \(A,C\) alone are unable to do so. 
The size of this gap is bounded from below
by the quantity \(2/\|H^{-1}\|\), cf.~(\ref{equality}).\\

As a preparation we will consider the matrices of the form
\(I + AC\) with \(A,C\) positive semidefinite. These will play a key role
in our estimates and may have an independent
interest of their own. 
Note that they are generally not Hermitian.
\begin{theorem}\label{th_I+AC}
Let \(A,C\) be Hermitian positive semidefinite. Then

\noindent(i)
\begin{equation}\label{sigmaAC}
\sigma(AC) = \sigma(CA) \subseteq [0,\infty),
\end{equation}
(ii)
\begin{equation}\label{weaker_bounds}
\|(I + AC)^{-1}\| \leq 1 +
\frac{\min\{\|A\|^{1/2}\|L^*C\|,\,\|C\|^{1/2}\|\|AM\|\}}
{1 + \min\sigma(AC)}
\end{equation}
where 
\begin{equation}\label{factors}
A = LL^*,\quad  C = MM^*,
\end{equation}
(iii) the matrix \(AC\) is diagonalisable.
\end{theorem}
{\bf Proof.} The statements (i), (iii) above are not new
(see \cite{hla_oml}, \cite{horn}, respectively).
To prove (ii)  note that
 \begin{eqnarray}\label{formula}
 (\lambda I + AC)^{-1}&=& \frac{1}{\lambda}I - 
\frac{1}{\lambda}(\lambda I + AC)^{-1}AC =
 \frac{1}{\lambda}I - \frac{1}{\lambda}(\lambda I + LL^*C)^{-1}LL^*C\nonumber\\
 &=&\frac{1}{\lambda}I - \frac{1}{\lambda}L(\lambda I + L^*CL)^{-1}L^*C,
\end{eqnarray}
So the spectra of \(CA\) and \(L^*CL\) coincide - up to possibly the point zero.
 Now, \(L^*CL\) is Hermitian positive semidefinite, hence 
  \[
  \|(I + L^*CL)^{-1}\| = 
  \frac{1}{1 + \min\sigma(L^*CL)}
  =  \frac{1}
  {1 + \min\sigma(AC)}
  \]
  and therefore
  \[
  \|(I + AC)^{-1}\| \leq
  1 + 
       \frac{\|L\|\|L^*C\|}
        {1 + \min\sigma(AC)}.
  \]
  The second half of (\ref{weaker_bounds}) is similar. Q.E.D.\\

  From (\ref{weaker_bounds}) it immediately follows that 

  \begin{equation}\label{invI+AC}
  \|(I + AC)^{-1}\| \leq 1 + 
\frac{\|A\|^{1/2}\|C\|^{1/2}\|L^*M\|}
  {1 + \min\sigma(AC)} \leq 1 + \|A\|\|C\|.
  \end{equation}
\begin{proposition}\label{pr_invI+AC}
Let \(A,C\) be Hermitian positive semidefinite. Then
\begin{equation}\label{norm_I+AC}
\|I + AC\| \geq 1
\end{equation}
and equality is attained if and only if \(AC = 0\).
\end{proposition}
{\bf Proof.} Since the norm dominates the spectral radius,
and by (\ref{sigmaAC}) the latter is not less than one,
(\ref{norm_I+AC}) follows. In the case of equality, 
the whole spectrum consists of the single point \(1\), 
that is, the spectrum of \(AC = A^{1/2}A^{1/2}C\) is \(\{0\}\). 
Then the spectrum of the Hermitian matrix \(A^{1/2}CA^{1/2}\) 
also equals \(\{0\}\). Hence \(A^{1/2}CA^{1/2} = 
(C^{1/2}A^{1/2})^* C^{1/2}A^{1/2} = 0\) and then also
\(AC = (CA)^* = 0\). Q.E.D.

\begin{theorem}\label{HBinv}
Let in (\ref{quasidef}) the matrix \(B\) be square and invertible
and let 
\begin{equation}\label{alpha_gamma}
\alpha = \sup_{x \neq 0}\frac{x^*Ax}{x^*\sqrt{BB^*}x},\quad
\gamma = \sup_{x \neq 0}\frac{x^*Cx}{x^*\sqrt{B^*B}x}.
\end{equation}
Then
\begin{equation}\label{eq:HBinv}
\|H^{-1}\| \leq 
\|B^{-1}\|\left(1 + \max\{\alpha,\gamma\} +
 \alpha\gamma\right).
\end{equation}
\end{theorem} 
{\bf Proof.} Using the polar decomposition \(B = U\sqrt{B^*B} = \sqrt{BB^*}U\)
we get the factorisation 
\begin{equation}\label{factor}
H =
\left[\begin{array}{rr}
(BB^*)^{1/4}  &       0 \\
0   &    (B^*B)^{1/4}   \\
\end{array}\right]
\left[\begin{array}{rr}
\tilde{A}     &       U \\
U^*   &      -\tilde{C}  \\
\end{array}\right]
\left[\begin{array}{rr}
(BB^*)^{1/4}  &       0 \\
0   &    (B^*B)^{1/4}   \\
\end{array}\right]
\end{equation}
with
\[
\tilde{A} = (BB^*)^{-1/4}A(BB^*)^{-1/4},\quad
\tilde{C} = (B^*B)^{-1/4}C(B^*B)^{-1/4}.
\]
Also,
\begin{equation}\label{inv_tilde}
\left[\begin{array}{rr}
\tilde{A}     &       U \\
U^*   &      -\tilde{C}  \\
\end{array}\right]^{-1}
=
\left[\begin{array}{rr}
(I + \hat{C}\tilde{A})^{-1} \hat{C}    &       U (I + \tilde{C}\hat{A})^{-1}\\
(I + \hat{A}\tilde{C})^{-1}U^*  &      -\hat{A}(I + \tilde{C}\hat{A})^{-1}  \\
\end{array}\right]
\end{equation}
where
\[
\hat{A} = U^*\tilde{A}U,\quad
\hat{C} = U\tilde{C}U^*
\]
are again Hermitian positive semidefinite. This is immediately
verified taking into account the identity

\begin{equation}\label{identity}
(I + \hat{A}\tilde{C})^{-1}U^* = U^*(I + \tilde{A}\hat{C})^{-1}.
\end{equation}
This, together with the identities of the type
\begin{equation}\label{identity2}
(I + \hat{C}\tilde{A})^{-1}\hat{C} =
\hat{C}^{1/2}(I + \hat{C}^{1/2}\tilde{A}\hat{C}^{1/2})^{-1}\hat{C}^{1/2}
\end{equation}
and the obvious inequality
\[
\left|\left|\left[
\begin{array}{ll}
E & F  \\    
G & K  \\
\end{array} 
\right]\right|\right|
\leq \max\{\|E\|,\|K\|\} + \max\{\|F\|,\|G\|\},
\]
permits the use of (\ref{invI+AC}) and the factorisation
(\ref{factor}) to obtain (\ref{eq:HBinv}).
Here we have used the obvious identities
\[
\alpha = \|\tilde{A}\| = \|\hat{A}\|,\quad
\gamma = \|\tilde{C}\| = \|\hat{C}\|
\]
and the fact that \(U\) is unitary. Q.E.D.\\

If \(B\) is replaced by \(tB\), \(t > 0\) then (\ref{eq:HBinv})
goes over into
\begin{equation}\label{eq:HtBinv}
\|H^{-1}\| \leq \frac{\|B^{-1}\|}{t}
\left(1 + \frac{\max\{\alpha,\gamma\}}{t} +
 \frac{\alpha\gamma}{t^2}\right).
\end{equation}
Note that here the right-hand side is monotonically decreasing in 
\(t\).

The proof of Theorem \ref{HBinv} may appear odd: the estimate
for the inverse of a Hermitian matrix \(H\) relies heavily on
the estimate for the inverse of some {\em non-Hermitian} matrices
of the type \(I + AC\). But this is the price for halving the dimension
of the problem in working with 'non-symmetric' Schur complements.

On the other hand, if both \(A\), \(C\) are positive definite
then setting \(\tilde{C} = \hat{C} = C,\
\tilde{A} = \hat{A} = A,\ U = I\) in (\ref{inv_tilde}), the inclusion 
(\ref{gapHH0}) yields
\begin{equation}\label{gapHH1}
\|(I + AC)^{-1}\| \leq \max\{\|A^{-1}\|,\|C^{-1}\|\}.
\end{equation}

\begin{rem}\rm
By (\ref{sigmaAC}) the spectrum of \(I + AC\) is uniformly bounded
away from zero, so one may ask whether there is a uniform bound
for the norm of its inverse. The answer is negative as shows the following 
example which is due to  M. Omladi\v{c} (private
communication). 
Set
\[
A = \left[
\begin{array}{ll}
t & 0  \\    
0 & 1/t  \\
\end{array} 
\right],\quad
C = \left[
\begin{array}{ll}
1/t & 1  \\
1 & t  \\
\end{array}
\right].  
\]
Then 
\begin{equation}\label{omla}
I + AC = 
\left[
\begin{array}{ll}
2   & t \\
1/t & 2 \\
\end{array}
\right],\quad
(I + AC)^{-1} = \frac{1}{3}
\left[
\begin{array}{rr}
2    & -t \\
-1/t &  2 \\
\end{array}
\right],
\end{equation}
and this is not bounded as \(t\) varies over the positive reals.
\end{rem}

Numerous numerical experiments with random matrices
led us to conjecture the bound
\begin{equation}\label{conjectured_bound}
\|(I + AC)^{-1}\| \leq \|I + AC\|.
\end{equation}
This conjecture is true (i) in dimension two, (ii)
if one of the matrices \(A,C\) has rank one and (iii)
if \(A,C\) commute; in the last case with the trivial bound
\begin{equation}\label{triv_AC}
\|\left(I + AC\right)^{-1}\| \leq 1.
\end{equation}
However, the estimate (\ref{conjectured_bound}) is in general false.
A nice counterexample, communicated to
us by A. B\"{o}ttcher, is as follows. Set
\[
\left[
 \begin{array}{rrr}
 1   &  0   &  0  \\
 -20 & 1.1  & 0   \\
0    & -20  & 1.2 \\
 \end{array}
 \right].
\]
A numerical calculation gives 
\[
\|I + M\| = 21.177 < 22,\quad \|(I + M)^{-1}\| = 43.774 > 42.
\]
Now, (cf.~eg.~\cite{ballantine}) any diagonalisable matrix \(M\) with
non-negative eigenvalues (our \(M\) is such) is a product of two Hermitian
positive semidefinite matrices. Indeed, if
\[
M = U\Lambda U^{-1}
\]
with \(\Lambda \geq 0\) diagonal then
\[
M = (U\Lambda^{1/2}U^*)(U^{-*}\Lambda^{1/2} U^{-1}),
\]
thus yielding a counterexample to the conjecture.\\

 We now turn to the more complicated case
in which \(A,C\) may have null-spaces and the invertibility of \(H\)
is due to the conspiring of all three blocks
\(A,B,C\). As an additional information we will need lower
bounds for the non-vanishing part of \(\sigma(A),\sigma(C)\).
 Thus, it will be technically convenient to represent \(H\) in a block
form which explicitly displays these null-spaces:
\begin{equation}\label{quasidef_blocks}
H =
\left[\begin{array}{rrrr}
A        &   0      &   B_{11} &  B_{12} \\
0        &   0      &   B_{21} &  B_{22} \\
B_{11}^* &  B_{21}^*&   -C     &   0     \\
B_{12}^* &  B_{22}^*&   0      &   0     \\
\end{array}\right].
\end{equation}
Here, for the notational simplicity, the new blocks \(A,C\)
are the `positive definite restrictions' of the original blocks \(A,C\)
in (\ref{quasidef}).
In view of Proposition \ref{singularity}, \(H\) is nonsingular
if and only if both matrices
\[
\left[\begin{array}{rr}
B_{21}^* &  B_{22}^* \\
\end{array}\right],\quad
\left[\begin{array}{r}
B_{12} \\
B_{22} \\
\end{array}\right]
\]
have full rank.
The following theorem gives a new sufficient condition for
invertibility and subsequently a gap estimate.
\begin{theorem}\label{zero_dichotomy} Suppose that 
\[
H =
\left[\begin{array}{rr}
A     &       B \\
B^*   &      -C  \\
\end{array}\right]
\]
is quasi-semidefinite. Assume, in addition,
\begin{enumerate}
\item \(\dim(\mathcal{N}(A)) =\dim(\mathcal{N}(C))\)
\item \(B\) is a one-to-one map from \(\mathcal{N}(C)\)
onto \(\mathcal{N}(A)\).
\end{enumerate}
That is, the block \(B_{22}\) in (\ref{quasidef_blocks})
is square and nonsingular.
Then
\[
(-\epsilon, \epsilon)\cap \sigma(H) = \emptyset
\]
with
\[
\epsilon =
\frac{1}
{(1 + \max\{\|B_{12}B_{22}^{-1}\|,\|B_{21}^*B_{22}^{-*}\|\})^2
\max\{\|A^{-1}\|,\|C^{-1}\|, \|B_{22}^{-1}\|\}}.
\]
\end{theorem}
{\bf Proof.} We represent \(H\) in the unitarily equivalent,
permuted form
\begin{equation}\label{quasidef_permuted}
\left[\begin{array}{llll}
A         &   B_{11} &  0       &  B_{12} \\
 B_{11}^* &  -C      &  B_{21}^*&    0    \\
0         &  B_{21}  &   0      &  B_{22}  \\
B_{12}^*  &    0     & B_{22}^* &   0      \\
\end{array}\right]
=
\left[\begin{array}{ll}
\hat{A}       &   \hat{B}\\
 \hat{B}^*    &   \hat{G}  \\
\end{array}\right].
\end{equation}
By renaming this matrix again into \(H\) we now perform the decomposition
\[
H =
\left[\begin{array}{ll}
I     &   \hat{B}\hat{G}^{-1}\\
0      &         I \\
\end{array}\right]
\left[\begin{array}{ll}
\hat{A}  - \hat{B}\hat{G}^{-1}\hat{B}^*  &     0     \\
0                                        &   \hat{G} \\
\end{array}\right]
\left[\begin{array}{ll}
I                      &   0  \\
\hat{G}^{-1}\hat{B}^*  &  I   \\
\end{array}\right].
\]
This yields the simple estimate
\[
\|H^{-1}\| \leq (1 + \|\hat{B}\hat{G}^{-1}\|)^2
\max\{\|(\hat{A}  - \hat{B}\hat{G}^{-1}\hat{B}^*)^{-1}\|,\|\hat{G}^{-1}\|\}.
\]
We now bound the single factors above:
\begin{eqnarray*}
\hat{A}  - \hat{B}\hat{G}^{-1}\hat{B}^* &=& 
\hat{A} - 
\left[\begin{array}{ll}
0     &   B_{12}  \\
B_{21}^*  &   0   \\
\end{array}\right]
\left[\begin{array}{ll}
0     &   B_{22}^{-*}  \\
B_{22}^{-1}  &   0   \\
\end{array}\right]
\left[\begin{array}{ll}
0     &   B_{21}  \\
B_{12}^*  &   0   \\
\end{array}\right]\\
&=&
\left[\begin{array}{ll}
A     &   B_{11} - B_{12}B_{22}^{-1}B_{21} \\
B_{11}^* - B_{21}^*B_{22}^{-*}B_{12}^*  &   -C   \\
\end{array}\right].
\end{eqnarray*}
This matrix is quasi-definite, hence
the interval \((-\min\sigma(C),\min\sigma(A))\)
contains none of its eigenvalues. Thus, 
\(\hat{A}  - \hat{B}\hat{G}^{-1}\hat{B}^*\) 
is invertible and
\[
\|(\hat{A}  - \hat{B}\hat{G}^{-1}\hat{B}^*)^{-1}\| \leq \max\{\|A^{-1}\|,\|C^{-1}\|\}.
\]
Furthermore,
\[
\|\hat{G}^{-1}\| = \|B_{22}^{-1}\|
\]
and
\[
\hat{B}\hat{G}^{-1} = 
\left[\begin{array}{ll}
B_{12}B_{22}^{-1}     &       0 \\
0    &    B_{21}^*B_{22}^{-*} \\
\end{array}\right],
\]
whence
\[
\|H^{-1}\| \leq 
(1 + \max\{\|B_{12}B_{22}^{-1}\|,\|B_{21}^*B_{22}^{-*}\|\})^2
\max\{\|A^{-1}\|,\|C^{-1}\|, \|B_{22}^{-1}\|\}.
\]
Q.E.D.\\

Note that the radius of the resolvent interval guaranteed
in the previous theorem depends on the spectra of
some operators obtained from the original blocks
\(A,B,C\). 

If in the preceding theorem we replace \(B\) by \(tB\) and \(t\)
is sufficiently large then we obtain
\begin{equation}\label{eq:zero_dichotomy}
\epsilon =
\frac{t}
{(1 + \max\{\|B_{12}B_{22}^{-1}\|,\|B_{21}^*B_{22}^{-*}\|\})^2
\|B_{22}^{-1}\|}.
\end{equation}
Another relevant special case has \(A = C\) and 
\(B^* = B\), {\em both} positive definite. Then, as was shown
in \cite{kirsch}, we have 
\begin{equation}\label{kirsch}
\rho(H) \supseteq
(-\sqrt{\min\sigma(A)^2 + \min\sigma(B)^2}, 
  \sqrt{\min\sigma(A)^2 + \min\sigma(B)^2}\ ).
\end{equation}
\begin{rem}\rm

The technique used in the proof of Theorem \ref{th_I+AC} is related to
the 
more general functional calculus for products \(AC\) with \(A,C\) bounded and
selfadjoint and \(C\) positive semidefinite in a general Hilbert space.
 It reads
\begin{equation}\label{f_calc}
f \mapsto f(AC) = f(0)I + AC^{1/2}f_1(C^{1/2}AC^{1/2})C^{1/2}
\end{equation}
with
\[
f_1(\lambda) = 
\left\{\begin{array}{r}
\frac{f(\lambda) - f(0)}{\lambda},\quad \lambda  \neq 0,\\
f'(0),\quad \lambda  = 0,\\
\end{array}
\right.
\]
By the property 
 \begin{equation}\label{f(XY)X}
f(XY)X = Xf(YX),
\end{equation}
valid for any matrix analytic function \(f\),
this obviously extends the standard analytic 
functional calculus and
requires \(f\) to be differentiable at zero and otherwise just to be 
bounded and measurable; then \(f_1\) will again be bounded and measurable
and is applied to a selfadjoint operator 
\(C^{1/2}AC^{1/2}\).\footnote{
    This functional calculus, probably well-known by now, was communicated
    to the second author by the late C. Apostol, Bucharest,
     more than forty years ago.} This calculus is a Hilbert-space
generalisation of the assertions of Theorem \ref{th_I+AC} (i), (ii),
only here the point zero may remain a sort of a `spectral singularity'.

The linearity and multiplicativity 
of the map \(f \mapsto f(AC)\) is
 verified by straightforward algebraic manipulations. Also, if the functions
\(f\) in (\ref{f_calc}) are endowed with the norm
\begin{equation}\label{f_norm}
\|f\| = |f(0)| + \|f_1\|_\infty
\end{equation}
 then the map \(f\mapsto f(AC)\) is obviously continuous.
    This admits estimating some other interesting functions of \(AC\),
    for instance, the group \(e^{-ACt}\), if both \(A\) and \(C\)
    are positive semidefinite and \(t > 0\). In this case 
    $f(\lambda)={\rm e}^{- t \lambda}$ and it is immediately 
    verified that \(|f_1|,\ \lambda \geq 0\), is bounded by \(t\), 
    whence  
\begin{equation}\label{expAC}
\|e^{-ACt}\| \leq 1 + t\|AC^{1/2}\|\|C^{1/2}\|
\end{equation}
and similarly
\begin{equation}\label{expCA}
\|e^{-ACt}\| \leq 1 + t\|CA^{1/2}\|\|A^{1/2}\|.
\end{equation}

\end{rem}
\begin{rem} \rm{\bf Extending monotonicity properties?} 
In the introduction we have 
stated two known monotonicity properties of the eigenvalues for some
affine quasi-definite pencils. It is natural
to try to extend this monotonicity to some neighbouring classes of matrix
families.
Some of our examples will be of the form
\begin{equation}\label{4x4}
H =
\left[\begin{array}{rr}
A  &  B  \\
B^*   & -A  \\
\end{array}\right]
\end{equation}
with \(2\times 2\)-matrices \(A^* = A\) and \(B^* = \pm B\) and no 
(semi)definiteness assumption whatsoever. Here
a straightforward calculation shows that the characteristic polynomial
is
\begin{equation}\label{4x4pol}
\lambda^4 - 2(\|A\|_F^2 + \|B\|_F^2)\lambda^2 + 
\left|\det(A - \sqrt{\mp 1}B)\right|^2
\end{equation}
where \(\|\cdot\|_F\) means the Frobenius or Hilbert-Schmidt norm.
This can be used to give a general formula for the roots explicitly,
see the Appendix. 

If in a quasi-definite matrix 
(\ref{quasidef}) the matrices \(A\) and \(C\) increase
(in the sense of forms), then the estimate (\ref{gapHH0}) 
certainly improves,
but does the gap at zero also necessarily increase? The answer is no
as the following example due to W.~Kirsch (private communication) shows. Set
\begin{equation}\label{myKirsch}
H =
\left[\begin{array}{rr}
A  &  B_t  \\
B_t   & -A  \\
\end{array}\right]
\end{equation}
with
\[
A =
\left[\begin{array}{rr}
2 &  -1 \\
-1  & 2\\
\end{array}\right],\quad
B_t =
\left[\begin{array}{rr}
1  &  0 \\
0  &  t \\
\end{array}\right].
\]
The matrix is quasi-definite. 
By (\ref{4x4pol}) the characteristic equation  
is readily found to be
\begin{equation}\label{characteristic}
\lambda^4 - (11 + t^2)\lambda^2 + 13 + 5t^2 + 2t = 0
\end{equation}
and the absolutely smallest 
eigenvalue is given in Figure \ref{fig_1} as function of \(t\),
\(5 < t < 20\),
(conveniently scaled) and it
shows a non-monotonic behaviour. Thus, there does not seem to
be a simple generalisation of the monotonicity property (A). On the
other hand, by the unitary similarity  
\begin{equation}\label{Kirsch}
\left[\begin{array}{rr}
A   &  B_t  \\
B_t  & -A  \\
\end{array}\right]
=
\frac{1}{2}
\left[\begin{array}{rr}
I   &   I  \\
I  &   -I  \\
\end{array}\right]
\left[\begin{array}{rr}
B_t  &  A  \\
A   & -B_t  \\
\end{array}\right]
\left[\begin{array}{rr}
I   &   I  \\
I  &   -I  \\
\end{array}\right]
\end{equation}
the same holds for the property (B).
\begin{figure}[htp]
\begin{center}
\includegraphics[width=8cm]{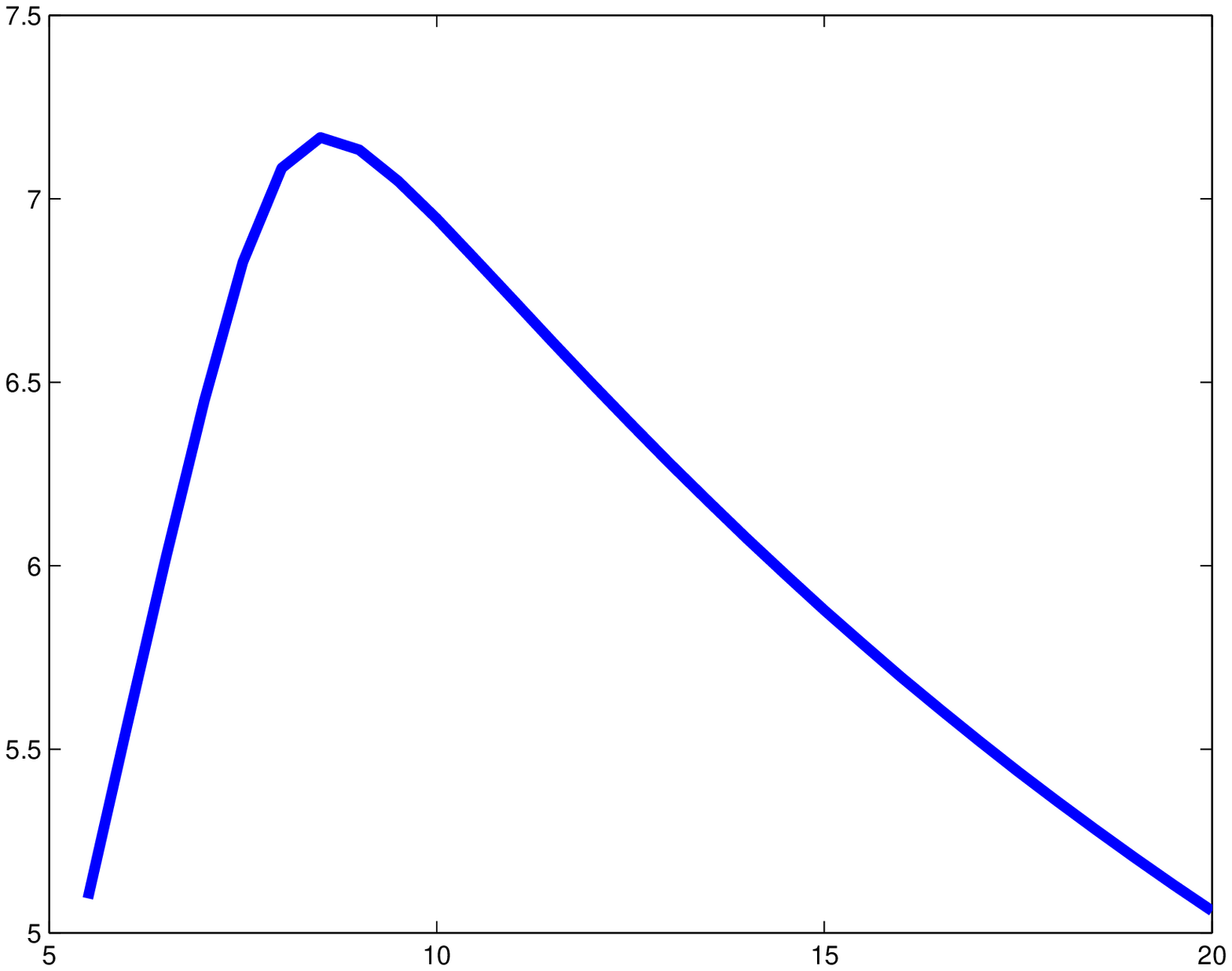}\caption{ \label{fig_1} 
Lack of monotonicity}
 \end{center}
 \end{figure}%
Another likely generalisation of (B), namely to have monotone 
eigenvalues if in 
(\ref{4x4}) the matrix
\(A\) is replaced by
\(tA\)
is also false. The counterexample is a numerical one:
\[
A =
\left[\begin{array}{rr}
1.24 &   0.81 \\
0.81 &   0.53\\
\end{array}\right],\quad
B =
\left[\begin{array}{rr}
0.30  &  -0.27 \\
-0.31 &  -0.48\\
\end{array}\right].
\]
Here both \(A\) and \(B\) are positive definite. 
The absolutely smallest eigenvalues  for \(5 < t < 20\)
are shown in the Figure
\ref{fig_2}.
\begin{figure}[htp]
\begin{center}
\includegraphics[width=9cm]{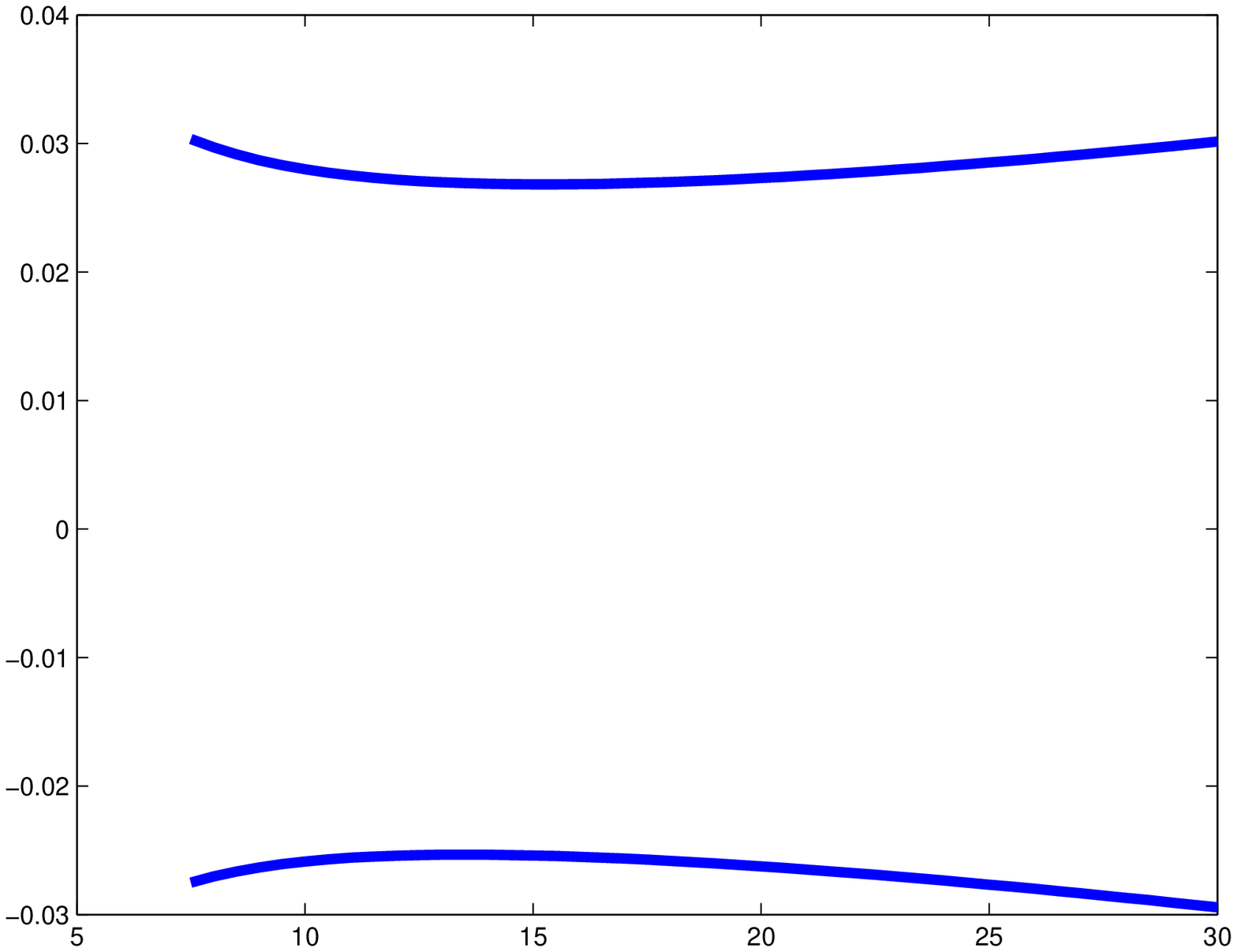}\caption{ \label{fig_2}
Lack of monotonicity}
 \end{center}
 \end{figure}%
\end{rem}
Finally, a yet simpler quasi-semidefinite example is given by 
\begin{equation}\label{simple}
H =
\left[\begin{array}{rr}
A_t  &  B  \\
-B   & -A_t  \\
\end{array}\right]
\end{equation}
with
\[
A_t =
\left[\begin{array}{rr}
0 &  0 \\
0  & t\\
\end{array}\right],\quad
B =
\left[\begin{array}{rr}
0  &  1 \\
-1  &  0 \\
\end{array}\right].
\]
Note that  \(\det(H) = 1\) for all \(t\) and that 
the spectrum $ -\lambda_1, -\lambda_2, \lambda_2, \lambda_1$ 
is symmetric w.r.t.~zero.
Thus if $|\lambda_1|$ increases, $|\lambda_2|$  has to decrease with growing $t>0$.
This already shows that property (B) in the introduction cannot hold 
for $t \mapsto A+t\left[\begin{array}{rr}
0 &  0 \\
0  & 1\\
\end{array}\right]$.
Moreover, it turns out that the spectral gap of \eqref{simple} shrinks to zero
as \(t \to \infty\).

At the end of this remark, let us formulate certain monotonicity properties of 
one-parameter families of quasi-definite matrices which are possibly true, 
but cannot prove at the moment.
The open questions are: 
\begin{itemize}
 \item If $A$ and $C$ in 
\begin{equation}\label{eq:Ht}
H_t =
\left[\begin{array}{rr}
A     &       tB \\
tB^*   &      -C  \\
\end{array}\right], \qquad t > 0
\end{equation}
are positive semidefinite, 
are all positive eigenvalues of $H_t$ isotone functions of $t$, 
and all negative eigenvalues of $H_t$ antitone functions of $t$? 
This would be an extension of property (A) mentioned in the introduction.
\item 
Are, in this situation, the positive eigenvalues of $H_t$ strictly increasing in $t$?
Under which conditions on $A$, $B$, and $C$?
\item 
Does this properties carry over to the infinite dimensional case, e.~g.~when 
$H_t $ in \eqref{eq:Ht} is defined on $\ell^2(\ZZ) \times \ell^2(\ZZ)$ 
and $A, B, C$ are bounded operators on $\ell^2(\ZZ)$?
\item 
A particularly interesting special class of operators of this type is
\begin{equation}\label{eq:Deltat}
H_t =
\left[\begin{array}{rr}
-\Delta&       tB \\
tB^*   &      -\Delta  \\
\end{array}\right], \qquad t > 0
\end{equation}
where $\Delta$ is the finite-difference Laplacian on $\ell^2(\ZZ)$, i.e.
\begin{equation*}
 \Delta \phi (x)= \sum_{y\in \ZZ, y\sim x} \big(\phi(y) -\phi(x)  \big ), \qquad x\in \ZZ, \quad \phi \in \ell^2(\ZZ)  
\end{equation*}
 and where $y\sim x$ denotes the neighbours of $x$.
\end{itemize}

\section{Unbounded operator matrices}\label{unbounded}
Most of the results obtained above immediately extend 
to infinite dimensional Hilbert space. 
Theorem  \ref{th_I+AC} (except (iii)) and Proposition \ref{pr_invI+AC}, 
together with their proofs, apply literally to any 
bounded selfadjoint positive semidefinite operators \(A,C\). 
Theorem \ref{HBinv} even allows \(B,A,C\) to be unbounded.
In fact, the last two may be just quadratic forms, requiring, 
of course, that the quantities \(\alpha,\gamma\) in (\ref{alpha_gamma}), 
reformulated in the quadratic form context, be finite, whereas \(B\) needs to have a bounded,
everywhere defined inverse.
More precisely, in this setting, the operator \(H\) is defined by the {\em form block matrix}
\begin{equation}\label{H_form}
\left[\begin{array}{rr}
\mathbf{a}    &       B \\
B^*   &      -\mathbf{c}  \\
\end{array}\right]
\end{equation}
where the symmetric sesquilinear forms \(\mathbf{a},\mathbf{c}\) have to be defined
on the form domains of \(\sqrt{BB^*},\sqrt{B^*B}\), respectively, and the
relative form bounds
\begin{equation}\label{alpha_gamma'}
\alpha = \sup_{x \neq 0}\frac{\mathbf{a}(x,x)}{\|(B^*B)^{1/4}x\|^2},\quad
\gamma = \sup_{x \neq 0}\frac{\mathbf{c}(x,x)}{\|(BB^*)^{1/4}x\|^2}
\end{equation}
need to be  finite.
 So, the operators
 \(\tilde{A},\tilde{C},\hat{A},\hat{C}\) appearing 
 in the proof of Theorem \ref{HBinv} will again be bounded and positive
semidefinite whereas the formula (\ref{factor}) now serves as a natural definition
of the operator \(H\) itself. Indeed, \(H\) is given  as a product
of three selfadjoint operators, each having a bounded, selfadjoint
inverse. The bounded invertibility of the first and the third factor in
(\ref{factor}) is trivial, whereas for the second it follows from the formula
(\ref{inv_tilde})  (cf. also \cite{ves}). Similar remarks hold for Theorem
\ref{zero_dichotomy} as well. (Proposition \ref{singularity} could also be 
reformulated in infinite dimension, but this will not interest us here.)

We will now compare our bound with a bound obtained in \cite{winklmeier}.
This bound (with our notation) requires
that \(A,C\) be relatively bounded with respect to \(B,B^*\),
respectively. According to \cite{kato}, Ch. VI, Th.~1.38, 
the operator boundedness implies the form boundedness with
the same bound; so our setting is more general. 
In addition  \cite{winklmeier} gives an eigenvalue bound under the
condition that at least one of the operators \(A,C\) is bounded.
The estimate obtained there is rather complicated; but if both \(A,C\)
are bounded then \cite{winklmeier} gives the somewhat simpler bound 
\begin{equation}\label{winkl_bd}
\dist(\sigma(H),0)\geq -\frac{1}{2}\left(\|A\| + \|C\|\right) + 
\left(\frac{1}{4}\left(\|A\| - \|C\| \right)^2 + \left(\|B^{-1}\|^{-2}
\right) \right)^{1/2},
\end{equation}
which is still not easily compared with our estimate (\ref{eq:HBinv}). 
Anyhow, if \(\|A\| = \|C\|\) and the relative bound 
\(\|A\| \|B^{-1}\|\) is larger than one, then the right-hand 
side of (\ref{winkl_bd})
becomes negative and the estimate is void
whereas the bound (\ref{eq:HBinv}) always makes sense.
In fact, our relative
bounds \(\alpha\) and \(\gamma\) may be arbitrary. In particular, they
need not to be less than one, which is a usual requirement in operator 
perturbation theory. It is a general feature with quasi-definite matrices 
that perturbations, as long as they respect, in an appropriate sense,
the block structure, need to be relatively bounded, but not necessarily
with the  relative bound less than  one, in order to yield an 
effective perturbation theory.
Such a phenomenon was already encountered in \cite{ves}, for example. 

The 
selfadjointness of the operator \(H\) from (\ref{H_form}) immediately applies to
various kinds of Dirac operators with supersymmetry 
(see \cite{thaller}, Sect.~5.4.2 and 5.5) under the appropriate definiteness
assumption for the diagonal blocks.

An analogous construction of a selfadjoint block operator matrix
\(H\) was made in \cite {ves} in the `dual' case in which
\(B\) is dominated by \(A,C\) in the sense that \(A^{-1/2}BC^{-1/2}\)
is bounded. 
Estimate (\ref{Hstretch}) extends to  this more general situation,
where $A,B,C$ need not be bounded.

Finally we come back to the estimate (\ref{kirsch}).
The proof given in \cite{kirsch} went through squaring the matrix 
\begin{equation}\label{Hequal}
H =
\left[\begin{array}{rr}
A   &        B \\
B   &       -A  \\
\end{array}\right],
\end{equation}
which is inconvenient if \(A,B\) are unbounded. We provide an alternate
proof under a weaker assumption, namely that instead of operators \(A,B\)
we have symmetric positive semidefinite (not necessarily closed) 
sesquilinear forms \(a,b\) defined on a dense domain 
\( {\cal D} ={\cal D}(a)={\cal D}(b)\).
The obvious generalisation of the block operator matrix
(\ref{Hequal}) is the symmetric sesquilinear form \(h\) defined
as
\begin{equation}\label{h}
h(x,y) = a(x_1,y_1) + b(x_2,y_1) + b(x_1,y_2) - a(x_2,y_2),\quad
x = 
\left[\begin{array}{rr}
x_1 \\
x_2 \\
\end{array}\right],\quad
y =
\left[\begin{array}{rr}
y_1 \\
y_2  \\
\end{array}\right]
\end{equation}
for
\[
x,y \in {\cal D}\oplus {\cal D}.
\]
Neither of the forms \(a,\ b\) need to be closed but 
{\em their sum} shall be assumed as closed.
Here we have, in fact, first to construct the operator \(H\).
To this end we use the `off-diagonalizing' transformation
given by the unitary matrix
\begin{equation}\label{U}
U =
\frac{1}{\sqrt{2}}
\left[\begin{array}{rr}
I    &        I \\
iI   &        -iI  \\
\end{array}\right]
\end{equation}
(cf.~\cite{thaller}).
Obviously \(U({\cal D}\oplus {\cal D}) 
= {\cal D}\oplus {\cal D}\), whereas a 
direct calculation leads to
\begin{eqnarray}
\hat{h}(x,y) &=& h(Ux,Uy) = ia(x_1,y_2) + b(x_1,y_2) + b(x_2,y_1) - ia(x_2,y_1)\nonumber\\
&=& \tau(x_2,y_1) + \tau^*(x_1,y_2)
\end{eqnarray}
where the forms
\begin{equation}\label{taus}
\tau = a -ib,\quad \tau^* = a +ib
\end{equation}
are sectorial and mutually adjoint. Obviously the range of the form
\(\tau\) lies in the lower right quadrant of the complex plane.

The form \(\tau\) is closed. This is readily seen from the 
equivalence of the corresponding norms:
\begin{eqnarray*}
|\tau(x,x)| &=& \sqrt{a(x,x)^2 + b(x,x)^2} \leq a(x,x) + b(x,x) \\
&\leq& \sqrt{2(a(x,x)^2 + b(x,x)^2)} 
=\sqrt{2}|\tau(x,x)|,
\end{eqnarray*}
so that the closedness of \(a + b\) is, in fact, equivalent to 
that of \(\tau\).
Thus \(\tau,\tau^*\) 
generate
mutually adjoint maximal sectorial operators
\(T,T^*\), respectively (see \cite{kato}, Ch.~VI, Theorem 2.1).
Now, for \(x_2 \in {\cal D}(T), x_1 \in {\cal D}(T^*)\) and
\(y_1,y_2 \in {\cal D}\) we have
\begin{equation}\label{hHhat}
\hat{h}(x,y) =  (\hat{H}x,y) 
\end{equation}
where the operator
\begin{equation}\label{Hhat}
\hat{H} =
\left[\begin{array}{rr}
0    &        T \\
T^*  &         0  \\
\end{array}\right]
\end{equation}
is obviously selfadjoint with the domain
\({\cal D}(T^*) \oplus {\cal D}(T)\). 
Also selfadjoint is its inverse conjugate
\[
H = U\hat{H}U^*
\]
with
\begin{equation}\label{H}
h(x,y) = (Hx,y),
\quad {\cal D}(H) \subseteq {\cal D}(\tau) \oplus {\cal D}(\tau).
\end{equation}
The operator  \(H\) is uniquely determined by (\ref{H})  
as is shown in \cite{ves}, Proposition 2.3.

To estimate the inverse note that
\[
\|Tz\|\|z\| \geq |(Tz,z)| = \sqrt{a(z,z)^2 + b(z,z)^2} \geq
\sqrt{\alpha^2 + \beta^2}\ \|z\|^2 
\]
where \(\alpha,\ \beta \geq 0\) is the lower bound of \(a,\ b\),
respectively. The above
`Lax-Milgram inequalities' are, in fact, the key argument
in this matter. They are non-trivial if any of \(\alpha,\ \beta \)
is different from zero. In this case, 
by the maximality of \(T\), its inverse is everywhere defined and
\[
\|T^{-1}\| = \|T^{-*}\| \leq \frac{1}{\sqrt{\alpha^2 + \beta^2}}, 
\]
From the obvious formula
\begin{equation}\label{Hhatinv}
\hat{H}^{-1} =
\left[\begin{array}{rr}
0    &        T^{-*} \\
T^{-1}  &         0  \\
\end{array}\right]
\end{equation}
we finally obtain
\begin{equation}\label{Hinvfinal}
\|\hat{H}^{-1}\| =  \|H^{-1}\|  = \|T^{-1}\| \leq 
 \frac{1}{\sqrt{\alpha^2 + \beta^2}}, 
\end{equation}
which obviously reduces to (\ref{kirsch}) if \(a,b\) are bounded. Thus,
we have proved the following theorem.
\begin{theorem}\label{th:kirsch} Let \(a,b\) be positive 
semidefinite symmetric
sesquilinear forms with the common dense domain \({\cal D}\) and
respective lower bounds \(\alpha,\beta\) and such that \(a + b\) 
is closed.
Then the form \(h\) from (\ref{h}) defines a unique selfadjoint operator \(H\)
with \({\cal D}(H) \subseteq {\cal D} \oplus {\cal D}\) and 
\(h(x,y) = (Hx,y)\)
for \(x \in {\cal D}(H),\ y \in {\cal D}\oplus {\cal D}\). 
Moreover, if any of
\(\alpha,\beta\) is non-zero then \(H\) has a bounded inverse with
\[
\|H^{-1}\|  \leq 
 \frac{1}{\sqrt{\alpha^2 + \beta^2}}.
\]
\end{theorem}
\begin{rem}\rm (i) The conditions of the preceding theorem are 
obviously fulfilled if one of the forms \(a,\ b\) is closed
and the other is relatively bounded with respect to the first.
Moreover, if, say, \(b\) is relatively bounded with respect to \(a\) then
\(b\) need not to be semidefinite; indeed the whole construction 
of \(H,\hat{H},T\) in the proof of the preceding theorem goes through
and we have
\begin{equation}\label{oneform}
\|\hat{H}^{-1}\| =  \|H^{-1}\|  = \|T^{-1}\| \leq 
 \frac{1}{\alpha},
\end{equation}
provided that \(a\) is positive definite.

(ii) The form \(\tau\) constructed in the proof of the preceding theorem
is not sectorial in the strict sense as defined in \cite{kato}
 because its range
does not lie symmetrically with respect to the positive real axis.
But, of course, the whole theory developed in \cite{kato} 
naturally and trivially extends to all
kinds of numerical ranges having semi-angle less than \(\pi/2\).
The standard form can be achieved simply by multiplying 
\(\tau\) with a phase factor
\[
e^{\frac{\pi}{4}i}\tau = \frac{1}{\sqrt{2}}
((a + b) + i(a - b)).
\]
The symmetric part of this form is closed, whereas the skew-symmetric
part is relatively bounded  with respect to the symmetric one, so it
is sectorial in the strict sense of \cite{kato}.

(iii) The obvious fact that the eigenvalues (whenever existing)
of \(H\) are \(\pm\) singular values of \(T\) may have advantage
in numerical computations with finite matrices. 
Firstly, the size of \(T\) is half
the size of \(H\) and, secondly, there is plenty of reliable
computational software to compute the singular values (and vectors)
of arbitrary matrices.

(iv) If \(a + b\) is only closable then its closure is again
of the form \(\tilde{a} + \tilde{b}\) where \(\tilde{a},\ \tilde{b}\)
are obtained by the usual limiting process and Theorem 
\ref{th:kirsch} applies. We omit the details.
\end{rem}

\section{Stokes matrices}\label{Stokes matrix}

If we set $C=0$ in \eqref{quasidef}, we obtain a \emph{Stokes matrix}.
Stokes matrices have been extensively studied, see \cite{RuWa}, \cite{Axel} and
the literature cited there. 
For $C=0$, we obviously have $AC=CA$. Consequently by (\ref{triv_AC}) the estimate
(\ref{eq:HBinv}) becomes 
\begin{equation}\label{eq:HBinv0}
\|H^{-1}\| \leq \|B^{-1}\|(1 + \alpha).
\end{equation}
A more careful inspection of formula (\ref{inv_tilde})
gives a tighter bound
\begin{equation}\label{eq:HBinv01}
\|H^{-1}\| \leq \|B^{-1}\|\frac{\alpha + \sqrt{\alpha^2 + 4}}{2}.
\end{equation}
In \cite{RuWa} the following spectral inclusion was proved. 

\begin{theorem}[\cite{RuWa}]
For positive definite \(A\) 
and \(B\) of full column rank, 
\begin{eqnarray}\label{eq:RuWa}
\sigma(H)  &\subseteq&  I_+ \cup I_-,\\
I_+  &=&  \left(\alpha_1,
\frac{\alpha_m + \sqrt{\alpha_m^2 + 4\beta_m^2}}{2}\right),\\
I_-  &=&  \left(\frac{\alpha_1 + \sqrt{\alpha_1^2 + 4\beta_m^2}}{2},
\frac{\alpha_m + \sqrt{\alpha_m^2 + 4\beta_1^2}}{2}\right),
\end{eqnarray}
where
$0 < \alpha_1\leq \cdots\leq\alpha_m
$ are the eigenvalues of \(A\) whereas
$0 < \beta_1\leq\dots\leq\beta_m
$ are the singular values of \(B\).
\end{theorem}
Under the same assumptions \cite{Axel} establishes the 
inclusion \eqref{eq:RuWa} with the intervals  \(I_\pm\) given by
\begin{eqnarray}\label{eq:Axel}
I_-  &=&  \left(\frac{-2\sigma_m\alpha_m}
{\alpha_1 + \sqrt{\alpha_1^2 + 4\sigma_m^2}},
\frac{\sigma_1\alpha_1}{\sigma_1 + \alpha_1}\right),\\
I_+  &=&  \left(\alpha_1,
\frac{\alpha_m + \sqrt{\alpha_m^2 + 4\sigma_m^2}}{2}\right),
\end{eqnarray}
where
$0 < \sigma_1\leq \cdots\leq\sigma_m
$ are the eigenvalues of \(B^*A^{-1}B\).

In the following we  partly improve and generalise the foregoing results.
For illustration purposes let us start with the \(2 \times 2\)-case, i.~e.
\begin{equation}\label{stokes22}
H =
\left[\begin{array}{rr}
a    &       b \\
\bar{b}   &  0 \\
\end{array}\right],\quad \min\{a,|b|\} > 0.
\end{equation}
The eigenvalues of $H$ are 
\begin{equation}\label{f2x2}
\lambda_\pm = f_\pm(a,t) =
\frac{a \pm \sqrt{a^2 + t^2}}{2},\quad \text{ where }t = 2|b|.
\end{equation}
The functions \(f_-,f_+\) have the following properties:
\begin{enumerate}
\item \(f_-(a,t) < 0 < f_+(a',t')\) for all $a,t,a',t'\in 
\mathbb{R}$, 
\item \(f_+\) is increasing in both variables \(a,t\),
\item \(f_-\) is increasing in \(a\) and decreasing in  \(t\). 
\end{enumerate}
Herewith a result for $n\times n$ matrices.
\begin{theorem}\label{th:stokes}
Let 
\begin{equation}\label{stokesH}
H =
\left[\begin{array}{rr}
A     &       B \\
B^*   &      0  \\
\end{array}\right]
\end{equation}
be an \(n\times n\) Hermitian matrix over the
field \(\mathbb{K} \in \{\mathbb{R}, \mathbb{C}\}\) 
such that  \(A\) is positive semidefinite of order \(m\) and
\begin{equation}\label{NAB}
\mathcal{N}(A)\cap \mathcal{N}(B^*) = \{0\}.
\end{equation}
Define $ p_+,p_-\colon \mathbb{S}^{m-1}\to \mathbb{R}$, 
where $\mathbb{S}^{m-1}$ is the unit sphere in $\mathbb{K}^m $, by
\begin{equation}\label{ppmDelta}
p_\pm(x) = \frac{x^*Ax \pm \sqrt{\Delta(x)}}{2},\quad 
\Delta(x) = (x^*Ax)^2 + 4x^*BB^*x,\quad \|x\|=1
\end{equation}
and
\[
p_+^+ = \max_{\|x\| = 1}p_+(x), \quad
p_+^- = \min_{\|x\| = 1}p_+(x), 
\]
\[
p_-^+ = \max_{\|x\| = 1}p_+(x), \quad
p_-^- = \min_{\|x\| = 1}p_-(x).
\]
Then the following hold.
\begin{description}
\item[Extremal eigenvalues:]
The points $p_+^+, p_+^- , p_-^+ , p_-^- $ are eigenvalues of $H$.
\item[Spectral inclusion:] 
\begin{equation}\label{signain+-}
\sigma(H) \subseteq I_-\cup I_+\cup \{0\} 
\end{equation}
where
\begin{equation}\label{eq:myI+-}
I_-  =  [p_-^-,p_-^+]  \text{ and } I_+  =  [p_+^-,p_+^+]
\end{equation}
and 
\begin{equation}\label{I+leqI-}
p_-^+ = \max I_-  \leq 0 \leq p_+^- =\min I_+,
\end{equation}
\begin{equation}\label{I+<I-}
p_-^+ < p_+^-.
\end{equation}
\item[Monotonicity:] 
Consider the eigenvalues of $H$ as functions of the  submatrices
\(A,B\). 
Then
all eigenvalues are non-decreasing with \(A\), whereas
the non-positive eigenvalues are non-increasing and the 
non-negative ones
non-decreasing with \(BB^*\).\footnote{The terms in(de)creasing for \(A\)
and \(B\) mean the quadratic forms \(x^*Ax\) and \(x^*BB^*x = \|B^*x\|^2\),
respectively.}
\end{description}
\end{theorem}
{\bf Proof.}
The eigenvalue equation for \(H\) is written as
\begin{eqnarray}
Ax + By &=& \lambda x,\\
B^*x   &=& \lambda y.
\end{eqnarray}
For \(\lambda  \neq 0\) these equations
are equivalent to 
\begin{equation}\label{quadratic}
(\lambda^2 I - \lambda A - BB^*)x = 0,\quad x \neq 0, \quad
y =B^*x/\lambda.
\end{equation}
By assumption (\ref{NAB}), for $x\neq 0$ we have \(\Delta(x) > 0\) 
and from (\ref{quadratic}) and \(\|x\| = 1\) it
follows
\begin{equation}\label{quadratic_eq}
\lambda^2x^*x   - \lambda x^*Ax - x^*BB^*x = 0.
\end{equation}
Therefore
\(\lambda \in \{p_+(x), p_-(x)\}\) and \(p_\pm\) are real-valued. 
Obviously
\begin{equation}\label{pf+-}
p_\pm(x) = f_\pm(x^*Ax,x^*BB^*x)
\end{equation}
with \(f_\pm\) from (\ref{f2x2}); then the properties 
(\ref{I+leqI-}), (\ref{I+<I-}) immediately follow from
the property 1 of the functions \(f_\pm\) from (\ref{f2x2}).
With the property \(\Delta(x) > 0\)
the matrix pencil \(\lambda^2 I - \lambda A - BB^*\) 
is called \emph{overdamped}. 
In \cite{duffin} a minimax theory for the eigenvalues of
overdamped pencils was established. According
to this theory 
there are minimax formulae for the eigenvalues
\[
\lambda_1^- \leq \cdots \leq \lambda_m^- \in I_-\quad
\lambda_m^+ \leq \cdots \leq \lambda_1^m \in I_+
\]
reading
\begin{equation}\label{maxmin}
\lambda_k^+ = \max_{S_k}\min_{x\in S_k\atop \|x\| = 1}p_+(x),
\end{equation}
\begin{equation}\label{minmax}
\lambda_k^- = \min_{S_k}\max_{x\in S_k\atop \|x\| = 1}p_-(x)
\end{equation}
where \(S_k\) varies over all \(k\)-dimensional subspaces of
 $\mathbb{K}^m$.
In particular,
\[
\lambda_1^+ = \max_{\|x\| = 1}p_+(x), \quad
\lambda_m^+ = \min_{\|x\| = 1}p_+(x), 
\]
\[
\lambda_1^- = \min_{\|x\| = 1}p_-(x), \quad
\lambda_m^- = \max_{\|x\| = 1}p_+(x).
\]
Thus, the boundary points of the two intervals \(I_-\) and  
\(I_+\) are eigenvalues, given by $p_+^+, p_+^-$, $p_-^+ , p_-^- $.
All other eigenvalues are in the specified range.
It remains to prove the monotonicity statement. 
It is an immediate consequence of the formulae (\ref{pf+-}),
(\ref{maxmin}) and \eqref{minmax} and the monotonicity properties of
the functions \(f_\pm\) in (\ref{f2x2}). Q.E.D.\\

By its very construction the interval \(I_+\) is minimal
among those which contain all positive eigenvalues $\lambda_m^+ \leq \cdots \leq \lambda_1^m$ of $H$. 
In an analogous sense  \(I_-\) is minimal, as well.
So they are included in those from (\ref{eq:RuWa}) as well as 
in those from (\ref{eq:Axel}).

On the other hand our intervals can be used as a 
source for new estimates.
Assume that \(B^*\) has full column rank (in which case \(\sqrt{BB^*}\)
is positive definite) and take \(\alpha\) as in (\ref{alpha_gamma}).
Then 
\begin{equation}\label{mynewStokes}
\left(\frac{-2\beta_1}{\alpha + \sqrt{\alpha^2 + 4}},
\beta_1 \right) \setminus \{0\}  \subseteq \rho(H)
\end{equation}
where \(\beta_1 = \min_{\|x\| = 1}\|B^*x\|\). 
Indeed, the inequality \(p_+ \geq \beta_1\) is trivial.
Using again the monotonicity properties of the function
\(f_-\) from (\ref{f2x2}) and taking \(\|x\| = 1\) we obtain
\begin{align*}
p_-(x)  &= \frac{x^*Ax  - \sqrt{(x^*Ax)^2 + 4x^*BB^*x}}{2}
\\
&\leq 
\frac{\alpha x^*(BB^*)^{1/2}x -  
\sqrt{(\alpha x^*(BB^*)^{1/2}x)^2 + 4x^*BB^*x}}{2}
\\
&\leq 
\frac{\alpha x^*(BB^*)^{1/2}x -  
\sqrt{(\alpha x^*(BB^*)^{1/2}x)^2 + 4(x^*(BB^*)^{1/2}x)^2}}{2}
\\
&= x^*(BB^*)^{1/2}x\frac{\alpha - \sqrt{\alpha^2 + 4}}{2}
\\
&\leq -\frac{2\beta_1}{\alpha + \sqrt{\alpha^2 + 4}},
\end{align*}
where we have first used (\ref{alpha_gamma}), then the obvious 
inequality
\[
x^*BB^*x \geq (x^*(BB^*)^{1/2}x)^2
\]
and finally the identity
\[
\min_{\|x\| = 1}x^*(BB^*)^{1/2}x = \min\sigma((BB^*)^{1/2}) =
\min(\sigma({BB^*}))^{1/2} = \min_{\|x\| = 1}\|B^*x\| = \beta_1.
\]
This proves (\ref{mynewStokes}).
Note that (\ref{eq:HBinv01}) exactly reproduces 
the lower edge of the spectral gap (\ref{mynewStokes})
while the upper edge of the gap is not described correctly by (\ref{eq:HBinv01}).

Another immediate consequence of the monotonicity properties of
the functionals \(p_\pm\) are perturbation bounds for the eigenvalues
of the perturbed matrix
\[
\hat{H} = H + \widetilde H
\]
with
\[
\widetilde  H =
\left[\begin{array}{rr}
\widetilde  A     &      \widetilde   B \\
\widetilde  B^*   &      0  \\
\end{array}\right]
\]
where
\[
|x^*\widetilde  Ax| \leq \eta x^*Ax, \quad \|\widetilde  Bx\| 
\leq \eta \|Bx\|,\quad \eta < 1.
\]
Then, as was shown in \cite{vesslap}, the eigenvalues
\(\hat{\lambda}_1^-, \ldots, \hat{\lambda}_m^-,\hat{\lambda}_1^+, \ldots, \hat{\lambda}_m^+\) 
of the perturbed  matrix $\hat{H}$ satisfy
\begin{eqnarray}\label{pertH+}
(1 - \eta)\lambda_k^+ \leq & \hat{\lambda}_k^+  &
\leq  (1 + \eta)\lambda_k^+,\\
\label{pertH-}
\frac{1 + \eta}{1 - \eta}\lambda_k^-  \leq & \hat{\lambda}_k^- &
\leq \frac{1 - \eta}{1 + \eta}\lambda_k^-.
\end{eqnarray}
\begin{remark} 
The interest in Stokes matrices stems form the fact 
that they are discrete analogs of Stokes operators.
A  Stokes operator has the form
\[
H_S = 
\left[\begin{array}{rr}
-\diverg a\grad     & -\grad  \\
 \diverg   &      0  \\
\end{array}\right].
\]
Here \(a\colon \Omega \to (0,\infty)\) is a positive function on some domain
\(\Omega \subseteq \mathbb{R}^n\) such that the inverse of \(-\diverg a\grad\)
in \(L^2(\Omega)\) is compact. 
Operator-theoretical facts about
Stokes operators are given e.g.~in \cite{Schmitz}.

Without having checked the details of proofs we intutively expect that for such operators the monotonicity as well as the continuity bounds (\ref{pertH+})
for the positive eigenvalues as functions of  \(a(\cdot)\) should hold as well. 
Thus, a perturbation \(\hat a(\cdot)=a(\cdot) + \widetilde 
a(\cdot)\) of \(a(\cdot)\) satisfying
\[
|\widetilde  a(x)| \leq \eta a(x),\quad \eta < 1
\]
would imply  (\ref{pertH+}).
\end{remark}


\section{Boundary conditions and invertibility - a case study.}
\label{boundary}
The most prominent example of an operator whose invertibility depends on boundary
conditions is the Laplacian on an interval with Dirichlet and Neumann
boundary conditions.
A deeper manifestation of this phenomenon is encountered in the spectral analysis
of Schr\"odinger operators with periodic potential. Such operators exhibit
a spectrum consisting of intervals, so-called
spectral bands.
If one restricts the Schr\"odinger operator originally defined on
$\mathbb{R}$ or $\mathbb{R}^d$ to a finite interval or cube, respectively,
it is desirable to preserve the periodic structure
of the original, unrestricted operator as much as possible. A restriction to
a finite cube with Dirichlet boundary conditions leads to spurious
eigenvalues located in the spectral gaps of the
original operator. A consistent way to avoid these boundary-induced
eigenvalues is to impose periodic or, more generally, quasi-periodic
boundary conditions. For such restrictions, the arising
spectrum is contained in the spectrum of the original operator; 
see \cite{Mezincescu} for an exposition for operators on $L^2(\mathbb{R})$.
In the context of periodic Schr\"odinger operators, 
spectral pollution in gaps is a well studied subject, see e.~g.~\cite{Cances}.

In this section we want to explore these ideas applied to a block-operator
investigated in the recent paper \cite{Stolz}. 
There the following matrix of order \(n = 2m\) is considered:
\begin{equation}\label{stolz}
H_c = H_c(n) =
\left[\begin{array}{rr}
A + 2cI     &       B \\
-B    &      -A - 2cI\\
\end{array}\right]
\end{equation}
with the \(m\times m\)-blocks
\begin{equation}\label{AB}
A =
\left[\begin{array}{ccccc}
0     &    1    &    0    & \dots  &  0       \\ 
1     &    0    &    1    & \dots  & \vdots   \\
0     &    1    & \ddots  &        & \vdots   \\
\vdots&         &         &    0   &   1      \\
0     & \dots   & \dots   &   1    &   0      \\
\end{array}\right],\quad
B =
\left[\begin{array}{rrrrr}
0     &    1    &    0    & \dots  &  0       \\ 
-1     &    0    &    1    & \dots  & \vdots   \\
0     &    -1    & \ddots  &        & \vdots   \\
\vdots&         &         &   0    &   1      \\
0     & \dots   & \dots   &   -1    &   0      \\
\end{array}\right]
\end{equation}
where \(c\) is any real number (the factor \(2\) is set by convenience)
and \(n = 2m\).

We will analyse the spectrum of $H_c$ and see that it exhibits
two spurious eigenvalues.  To remove these, we will introduce a low-rank modification $\tilde H_c$ 
concentrated on the ``boundary''. This results in a circulant-type matrix. The circulant structure 
can be understood as an analogy to periodic boundary conditions used in the context 
of periodic Schr\"odinger operators. The specific type of the circulant matrix
shows that the operator considered in \cite{Stolz} lives on the two-fold covering space
$\{1,\ldots,2m\} \to \{1,\ldots,m\}$.

Of course, the spectrum of \(H_c = H_c(n) \) will depend on the 
dimension \(n\), so we will say
that an interval \({\cal I}\) around zero is a (maximal) {\em stable spectral gap} of \(H_c\)
if \( {\cal I}\cap \sigma(H_c) = \emptyset\) for all \(n\) and  \({\cal I}\) is
the largest interval with this property.

We perform the off-diagonalisation by taking the unitary matrix
\begin{equation}\label{UU}
U =
\frac{1}{\sqrt{2}}
\left[\begin{array}{rr}
I     &       I \\
I     &      -I \\
\end{array}\right]
\end{equation}
and obtaining
\begin{equation}\label{Kc}
K_c = U^{-1}H_cU = UH_cU =
\left[\begin{array}{rr}
0    &  A_c - B \\
A_c + B   &  0 \\
\end{array}\right]
=
2
\left[\begin{array}{rr}
0    &  T_c \\
T_c^*  &  0 \\
\end{array}\right]
\end{equation}
with
\begin{equation}\label{Tc}
T_c =
\left[\begin{array}{ccccc}
c    &        &           &     &      \\ 
1    &   c    &           &     &      \\
     &    1   &\ddots     &     &       \\
     &        & \ddots    &  c  &       \\
     &        &     &  1  &  c   \\
\end{array}\right]
\end{equation}
(all void places are zeros). As is well known, the eigenvalues of
\(K_c\), including multiplicities, are \(\pm\) the singular values of
\(T_c\). Now, the latter are of some importance in Matrix Numerical
Analysis, see \cite{erxiong}, where it was shown that the smallest
singular value of \(T_c\) tends to zero for \(m \to \infty\) and 
any fixed \(c\) with 
\(|c| < 1\). In any case the singular values of \(T_c\) are
independent of the sign of \(c\) as is seen from the property
\[
U_0T_cU_0 = T_{-c},\quad  [U_0]_{ij} = (-1)^j\delta_{ij}.
\]
We will now study these singular values in some detail. We shall distinguish
the cases

\begin{center}
(i) \(c = 0\), (ii) \(0 < c < 1\), (iii) \(c = 1\), (iv) \(c > 1\).
\end{center}

For \(c = 0\), the matrices \(T_c,T_c^*\) are partial isometries with all singular values 
equal to \(1\), except for the non-degenerate eigenvalue zero corresponding to
\[
{\cal N}(T_c) = \text{span } \{e_m\},\quad {\cal N}(T_c^*) = \text{span } \{e_1\}
\]
where \(e_1,\ldots,e_m\) is the canonical basis in \(\mathbb{R}^m\).
Hence \(K_c\) has the eigenvalues \(\pm 2\) each with multiplicity
\(n - 1\) and the double eigenvalue zero with
\begin{equation}\label{span2}
{\cal N}(K_c) = 
\text{span } 
\left\{\left[\begin{array}{r}
0\\
e_m   \\ 
\end{array}\right]
\left[\begin{array}{r}
e_1\\
0   \\ 
\end{array}\right]
\right\}.
\end{equation}
The case \(c > 1\) is easily accessible based on the representation
(\ref{stolz}), because then the matrix
\[
A_c := A + 2cI =
\left[\begin{array}{ccccc}
2    &    1    &    0    & \dots  &  0       \\ 
1     &   2    &    1    & \dots  & \vdots   \\
0     &    1    & \ddots  &        & \vdots   \\
\vdots&         &         &   2   &   1      \\
0     & \dots   & \dots   &   1    &  2      \\
\end{array}\right]
+ (2c -2)I
\]
is positive definite, being a sum of two obviously positive 
definite matrices, so \(\sigma(A_c) \geq 2c -2\). Therefore
\(H_c\) is quasidefinite and by 
(\ref{gapHH0}) the interval
\[
( -2c + 2,  2c -2)
\]
is contained in the stable spectral gap of \(H_c\) {\rm(}and of \(K_c\){\rm)}.
For further investigation we use the
fact that the singular values of \(T_c\) are the  square roots of the 
eigenvalues of \(T_c^*T_c\) or, equivalently, of
\begin{equation}\label{Wc}
W_c = T_{-c}^*T_{-c} =
\left[\begin{array}{rrrrr}
c^2 +1  &   -c     &    0    & \dots  &  0       \\ 
-c      &  c^2 +1  &  -c     & \dots  & \vdots   \\
0       &    -c    & \ddots  &        & \vdots   \\
\vdots  &          &         & c^2 +1 &   -c      \\
0       &  \dots   & \dots   &   -c   &  c^2       \\
\end{array}\right].
\end{equation}
Now \(W_cx = \lambda x\) is componentwise written as
\begin{eqnarray*}
(c^2 + 1)x_1 - cx_2 & = & \lambda x_1\\
-cx_{j-1} + (c^2 + 1)x_j - cx_{j+1} & = & \lambda x_j,\quad 
j = 2,\ldots, m-1,\\
-cx_{m - 1} + c^2x_m & = & \lambda x_m 
\end{eqnarray*}
or as a standard second order difference equation
\begin{equation}\label{diff}
-cx_{j-1} + (c^2 + 1 - \lambda)x_j - cx_{j+1} = 0, \quad j = 1,\ldots, m
\end{equation}
with the boundary conditions
\begin{equation}\label{bd}
x_0 = 0,\quad x_m - cx_{m+1} = 0.
\end{equation}
The solutions 
\begin{equation}\label{sinaj}
x_j = \sin j\alpha, \mbox{ with } \lambda = c^2 + 1 -2c\cos\alpha 
\end{equation}
and
\begin{equation}\label{sinhaj}
x_j = \sinh j\alpha, \mbox{ with } \lambda = c^2 + 1 -2c\cosh\alpha
\end{equation}
automatically  satisfy \(x_0 = 0\).
The second boundary condition from (\ref{bd}) will determine the
values of \(\alpha\). This gives
\begin{equation}\label{sinbd}
(1 - c\cos\alpha)\sin m\alpha - c\cos m\alpha \sin\alpha = 0,
\end{equation}
\begin{equation}\label{sinhbd}
(1 - c\cosh\alpha)\sinh m\alpha - c\cosh m\alpha \sinh\alpha = 0,
\end{equation}
respectively. 
In the easiest case \(c = 1 \), the substitution \eqref{sinaj} 
immediately leads to
\begin{equation*}
\alpha = \alpha_k = \frac{2k - 1}{2m + 1}\pi,\quad k = 1,\ldots,m,
\end{equation*}
giving rise to the eigenvalues
\begin{equation*}
\lambda = \lambda_k = 4 \sin^2\frac{2k - 1}{2m + 1}\pi ,\quad k = 1,\ldots,m. 
\end{equation*}
In this case the lowest eigenvalue
$\lambda_1 = 4 \sin^2\frac{2}{2m + 1}\pi \approx (2m + 1)^{-2}$
tends to zero as \(m\) tends to infinity,
so the stable spectral gap of $H_c$ is empty.

In the case \(c > 1\), equation (\ref{sinbd}) can be written as
\begin{equation}\label{tan}
f(\alpha) = \tan m\alpha\frac{1 - c\cos\alpha}{\sin\alpha} - c = 0,
\quad 0 < \alpha < \pi.
\end{equation}

\begin{figure}[htp]
\begin{center}
\includegraphics[width=11cm]{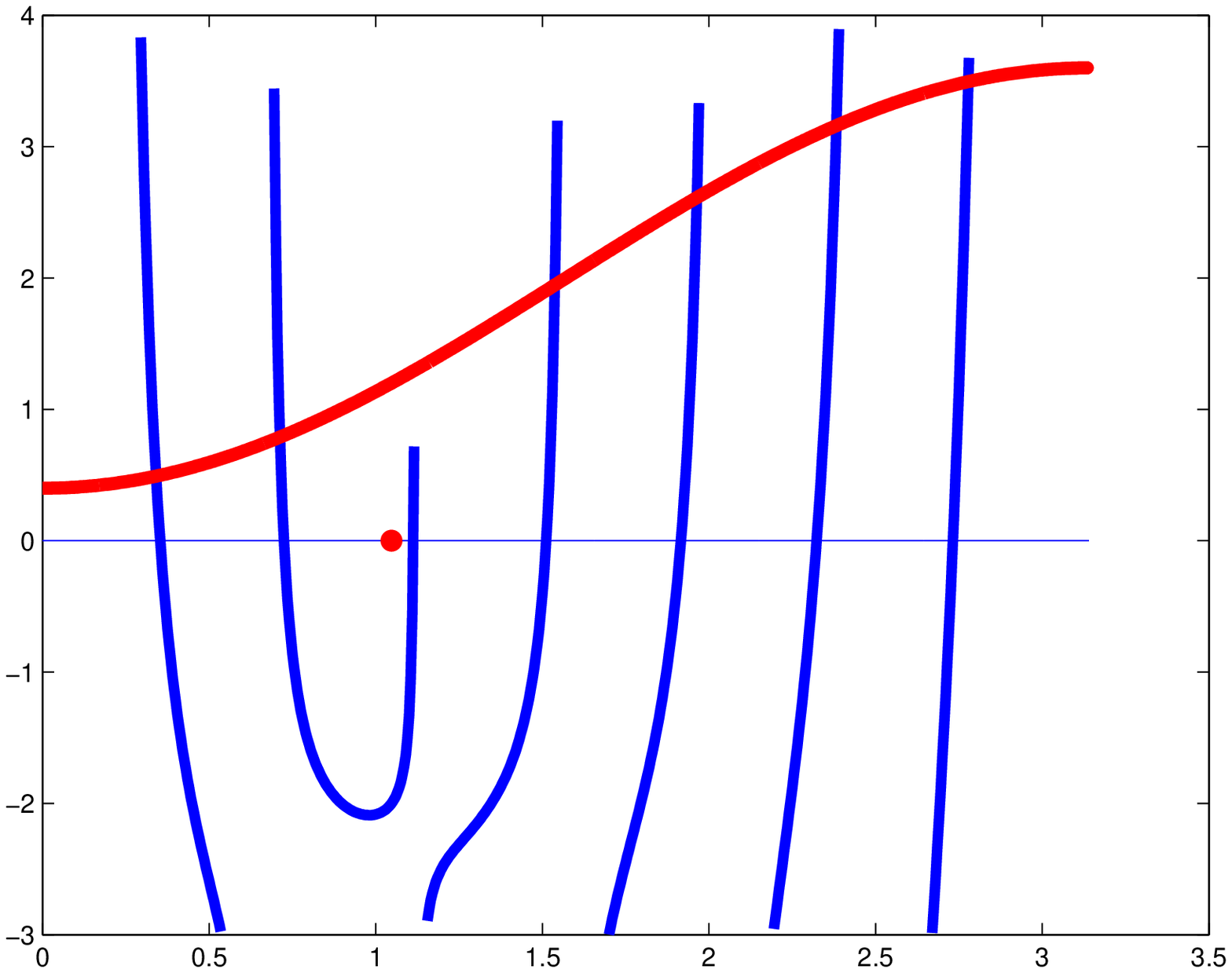}\caption{ \label{poles} 
Functions \(f\) and \(\lambda\) for \(n = 7\) and \(c = 2\)}
 \end{center}
 \end{figure}

The localisation of these roots is a bit involved, 
because there are several different cases to be distinguished. A
generic situation is shown on Figure \ref{poles}, which
displays
\begin{itemize}
\item the function \(f(\alpha)\) (blue) with its poles and roots,
\item the function \(\lambda =  \lambda = c^2 + 1 -2c\cos\alpha\)
(red), which, taken at the roots, gives the eigenvalues,
\item the point \(\hat{\alpha} = \arccos\frac{1}{c}\) on the 
 \(\alpha\)-axis on which \(f\) is generically negative.
\end{itemize}

Thus, between each two poles there is a root, 
except for the two poles enclosing \(\hat{\alpha}\);
these two poles enclose two roots, altogether \(n\) of them.
The case in which \(\hat{\alpha}\) coincides with one of the poles must be 
treated separately.
But to determine the exact position
of the stable spectral gap, it is enough to notice
that in any case, for \(m\) large enough, the interval
\[
0 < \alpha < \hat{\alpha},
\]
on which the factor \(1 - c\cos\alpha\) is positive,  
will contain several of the singularities 
\[
\frac{(2k -1)\pi}{2m},\quad k = 1,2,\ldots m
\]
of the function \(f\) in (\ref{tan}), see Figure \ref{poles}.
Each two of these singularities enclose  a root \(\alpha\)
of (\ref{tan}), and since each of them tends to zero for $m\to\infty,$
the lowest root tends to zero as well. Hence the corresponding
eigenvalue \(\lambda\) from (\ref{sinaj}) tends to \((c - 1)^2\).
Since we already know that the interval
\((-2(c - 1), 2(c - 1))\) is contained in the stable spectral gap of \(K_c\) (and \(H_c\)),
this interval is, in fact,  \emph{equal} to this gap.

In the case \(c < 1 \), the factor  \(1 - c\cos\alpha\) in (\ref{tan})
is globally positive, so the \(m\) singularities
\[
\frac{(2k -1)\pi}{2m}, \quad k = 1,2,\ldots, m
\]
will enclose \(m -1\) roots of the equation (\ref{tan}) and the lowest
of them will again approach zero with growing \(m\). That is,
the corresponding \(m -1\) eigenvalues \(\lambda\) from (\ref{sinaj}) will be 
larger than \((c - 1)^2\) and the lowest of them will approach \((c - 1)^2\).
Thus, we would have a stable spectral gap \((-2(1 - c), 2(1 - c))\),
but for one eigenvalue which is still to be determined. However,
as we know from \cite{erxiong}, the smallest eigenvalue tends
to zero with growing \(m\). This completes the picture. Thus, {\em for
\(0\leq c < \infty\) the interval
\[
(-2|c - 1|, 2|c - 1|)
\]
is the stable spectral gap for \(K_c\) {\rm(}and \(H_c\){\rm)}, 
{\bf except}
that for \(0 < c < 1\) there are two `spurious eigenvalues' tending
to zero with growing \(m\)}.

There is some interest in obtaining an asymptotic estimate
of the small eigenvalues. In \cite{erxiong} it was shown that
the smallest singular value of \(T_c\) is bounded from above by
\({\cal O}(\frac{1}{m})\). Numerical experiments indicate that
the decay is, in fact, much faster. The
setting of difference equations makes it possible to determine
the decay  accurately, and this is what we shall do now.

This solution is obtained by the ansatz \(x_j = \sinh j\alpha\)
and the fact that (\ref{sinhbd}) can be written as
\begin{equation}\label{tanh}
g(\alpha) =
\tanh m\alpha\frac{1 - c\cosh\alpha}{\sinh\alpha} - c = 0,
\quad 0 < \alpha < \infty.
\end{equation}
Since
\[
g(0_+) = m(1 - c) - c,\quad g(\infty) = - 2c
\]
and \(0 < c < 1\), for large \(m\) the equation (\ref{tanh}) has a 
positive root \(\alpha = \alpha_1\) which completes the \(m-1\) 
roots previously found, whereas the corresponding eigenvalue
\(\lambda = \lambda_1\) is given by (\ref{sinhaj}). 
As an approximation to  \(\alpha_1\) we propose the value
\begin{equation}\label{eq:alpha0}
\alpha_0 = \arcosh\frac{c^2 + 1}{2c}.
\end{equation}
Then a straightforward calculation gives 
\[
g(\alpha_0) = \frac{-2ce^{-2m\alpha_0}}{1 + e^{-2m\alpha_0}} =
-2ce^{-2m\alpha_0} + {\cal O}(e^{-2m\alpha_0})
\]
and
\[
g'(\alpha_0) = \frac{mc}{\cosh^2 m\alpha_0} 
- \tanh m\alpha_0\frac{-2c}{1 - c^2} =
-\frac{2c}{1 - c^2} + {\cal O}(e^{-2m\alpha_0})
\]
Thus, the difference \(\delta\alpha = \alpha_1 - \alpha_0\)
is given by
\[
\delta\alpha = - \frac{g(\alpha_0)}{g'(\alpha_0)} + {\cal O}(e^{-2m\alpha_0})
= {\cal O}(e^{-2m\alpha_0})
\]
whereas the corresponding eigenvalue \(\lambda_1\)  of $W_c=W_c(m)$ is given by
\[
{\cal O}(e^{-2m\alpha_0}) + \lambda_1  =
\delta\alpha(-2\sinh\alpha_0) =
- \frac{4c\delta\alpha}{1 - c^2} = 4ce^{-2m\alpha_0},
\]
where we have used the fact that the function \(c^2 + 1 -2c\cosh\alpha\)
vanishes at \(\alpha = \alpha_0 \).
Hence the small eigenvalues of \(K_c\) (and \(H_c\)) are asymptotically
absolutely bounded by
\[
{\cal O}(e^{-m\alpha_0}).
\]

We summarize the main findings in the following. 
\begin{theorem}
The spectrum of $H_c$ is symmetric w.r.t.~zero, i.e.~if $\lambda$ is an eigenvalue, then 
$-\lambda$  is also an eigenvalue, with the same multiplicity.

If $c\geq1$, the interval $(-2|c - 1|, 2|c - 1|)$ is a stable spectral gap,
i.e.~$(-2|c - 1|, 2|c - 1|)\cap \sigma(H_c(m))=\emptyset$ for all $m\in\NN$.

For $c\in[0,1]$ and each $m\in\NN$, 
$(-2|c - 1|, 2|c - 1|)\cap \sigma(H_c(m))$ consists of exactly two eigenvalues
with  absolute value of order ${\cal O}(e^{-m\alpha_0})$, with $\alpha_0$ as in \eqref{eq:alpha0}.

In both cases, $(-2|c - 1|, 2|c - 1|)$ is the  
\emph{maximal} interval with the above properties. 
More precisely, for any $\epsilon >0$ and  $N\in\NN$, there exists an $M\in\NN$ such that 
$(-2|c - 1|-\epsilon, 2|c - 1|+\epsilon)$ contains $2N$ eigenvalues  
of $H_c(m)$ for all $m\geq M$.
\end{theorem}

The spurious eigenvalues can be computed with high relative accuracy
by iteratively solving the equation (\ref{tanh}). By {\rm high relative
accuracy} we mean to obtain a significant number of correct digits
{\em independently} of the size of the computed quantity. Note that
the usual matrix computing software computes a
singular value of a matrix \(A\) with the error of the order 
\(\epsilon\|A\|\) (\(\epsilon\)
the machine precision)
which may yield no
significant digits, if the singular value itself is very small. 
A notable exception are bidiagonal matrices, which
is the case with our  \(T_c\). Then there exists an algorithm (and it is 
implemented in LAPACK and MATLAB packages) which computes each singular value
with about the same number of significant digits, no matter how small 
or how large it may be (barring underflow).\footnote{In fact, in order to
perform the computation with high relative accuracy,
MATLAB will need the input matrix to be {\em upper} bidiagonal, so the
MATLAB function {\em svd} should be applied not to
\(T_c\) but to its transpose.} 

It is also worthwhile to note that the components \(x_j = \sinh j\alpha_1\)
of the corresponding eigenvector
always agglomerate on one side of the sequence \(1,\ldots,m\), that is,
on the boundary, while all other eigenvectors exhibit standard oscillatory
behaviour. \\

{\bf Removing spurious eigenvalues.} The form of the 
null-space of \(K_0\) suggests to introduce the matrix
\begin{equation}\label{KK}
\tilde{K}_c = K_c + 
\left[\begin{array}{rr}
e_m    &       0 \\
0    &      e_1 \\
\end{array}\right]
\left[\begin{array}{rr}
-2     &       0 \\
0    &      2 \\
\end{array}\right]
\left[\begin{array}{rr}
e_m^T    &       0 \\
0    &      e_1^T \\
\end{array}\right]
= 
2
\left[\begin{array}{rr}
e_me_m^T    &      T_c\\
T_c^*   &      -e_1e_1^T \\
\end{array}\right].
\end{equation}
For \(c = 0\) this leaves all eigenvectors of \(K_0\) unchanged and raises
the zero eigenvalues to \(\pm 2\), respectively, thus `purging' the spurious
eigenvalues. The spectrum of \(\tilde{K}_0\) is \(\{\pm 2\}\) with the multiplicity \(m\). In particular,
\begin{equation}\label{KK2}
\tilde{K}_0^2 = 4I.
\end{equation}
It is a remarkable fact that for \(c \neq 0\) {\em the eigenvalues 
of the matrix \(\tilde{K}_c\) still come in
plus/minus pairs, including multiplicity}. To see this we take
\begin{equation}\label{J}
J = \left[\begin{array}{cccc}
0       & \cdots &  0     &  1    \\
0       & \cdots &  1     &  0    \\
\vdots  &        &        & \vdots\\
1       &   0    & \cdots &  0    \\
\end{array}\right]
\end{equation}
and set
\[
\mathbb{P} = 
 \left[\begin{array}{cc}
I       & 0    \\
0      &  J   \\
\end{array}\right]
\]
and note that \(J e_1 = e_m\), \(J e_m = e_1\) and that the matrix
\(S_c = T_cJ\) is symmetric.
Then
\[
\mathbb{P}^{-1}\tilde{K}_c\mathbb{P} = \mathbb{P}\tilde{K}_c\mathbb{P} =
 \left[\begin{array}{cc}
e_1e_1^T & S_c   \\
S_c        & -e_1e_1^T   \\
\end{array}\right].
\]
This matrix is of the form  (\ref{Hequal}), and such matrices
have the eigenvalues in plus/minus pairs when \(A\) and \(B\) are 
allowed to be any symmetric matrices (\cite{kirsch}).
It remains to determine the stable spectral gap of \(\tilde{K}_c\).
In order to do this it is convenient to turn back to the original
representation (\ref{stolz}) and to form the matrix
\begin{equation}\label{HH}
\tilde{H}_c = U\tilde{K}_cU =
\left[\begin{array}{rr}
\tilde{A} + 2cI  &  \tilde{B}  \\
\tilde{B}^*   &  \hat{A} - 2cI \\
\end{array}\right]
\end{equation}
with
\[
\tilde{A} =
\left[\begin{array}{cccccc}
1     &    1    &    0    &        &        &          \\ 
1     &    0    &    1    &        &        &          \\
      &    1    & \ddots  & \ddots &        &          \\
      &         & \ddots  & \ddots & \ddots &          \\
      &         &         & \ddots &   0    &    1      \\
      &         &         &        &   1    &   -1      \\
\end{array}\right],\quad
\tilde{B} =
\left[\begin{array}{cccccc}
1     &    1    &    0    &        &        &          \\ 
-1    &    0    &    1    &        &        &          \\
      &    -1   & \ddots  & \ddots &        &          \\
      &         & \ddots  & \ddots & \ddots &          \\
      &         &         & \ddots &   0    &   1      \\
      &         &         &        &  -1    &   
      1      \\
\end{array}\right],\quad
\]
\[
\hat{A} =
\left[\begin{array}{cccccc}
1     &   -1    &     0   &        &        &          \\ 
-1    &    0    &    -1   &        &        &          \\
      &    -1   & \ddots  & \ddots &        &          \\
      &         & \ddots  & \ddots & \ddots &          \\
      &         &         & \ddots &   0    &  - 1      \\
      &         &         &        &   -1   &   -1      \\
\end{array}\right].
\]
By (\ref{KK2}) we have \(\tilde{H}_0^2 = 4I\),
 which implies
\[
\tilde{A}^2 +  \tilde{B} \tilde{B}^* = 4I, \quad
\tilde{A}\tilde{B} +  \tilde{B}\hat{A} = 0, \quad
 \tilde{B}^*\tilde{B} + \hat{A}^2 = 4I.
\]
Using this\footnote{Of course, these three identities could be
proved directly.} we obtain that
\begin{equation}\label{HH2}
\tilde{H}_c^2 =
\left[\begin{array}{rr}
(4 + 4c^2)I + 4c\tilde{A}       &   0  \\
0        & (4 + 4c^2)I - 4c\hat{A}   \\
\end{array}\right].
\end{equation} 
Noting the identity
\[
J\hat{A}J = - \tilde{A}
\]
we obtain that the eigenvalues of \(\tilde{H}_c^2\) are 
\[
4 + 4c^2 - 4c\kappa_j,\quad j = 1,\ldots,r
\]
(with the multiplicity two) where \(\kappa_j\) are the eigenvalues
of \(\hat{A}\). These are obtained from the difference equation
\begin{equation}\label{diffAhat}
x_{j-1} + \kappa x_j + x_{j+1} = 0,\quad j = 1,\ldots,m
\end{equation}
with the boundary conditions
\begin{equation}\label{bdAhat}
x_0 = -x_1,\quad x_{m+1} = x_m.
\end{equation}
The substitution
\[
x_j = A\cos j\alpha + B\sin j\alpha 
\]
solves (\ref{diffAhat}) with
\[
\kappa = -2\cos\alpha,
\]
whereas the boundary conditions (\ref{bdAhat}) yield after
some computation
\[
\alpha = \alpha_k = \frac{2k -1}{2m}\pi, \quad k = 1,\ldots,m
\]
and hence (cf.~the Appendix)
\begin{equation}\label{kappa}
\kappa = \kappa_k = -2\cos\frac{2k -1}{2m}\pi, \quad k = 1,\ldots,m.
\end{equation}
Thus the eigenvalues of \(\tilde{H}_c^2\) are
\[
4 + 4c^2 + 2c\cos\frac{2k -1}{2m}\pi, \quad k = 1,\ldots,m
\]
(each taken twice), which is always larger than \(4(1 - c)^2\), and for
\(m\) large the set of these eigenvalues comes arbitrarily close to
\(4(1 - c)^2\). Again we conclude the following

\begin{theorem}
The stable spectral gap of \(\tilde{H}_c\) is
\[
(-2|c - 1|,2|c - 1|).
\]
More precisely, for all $m \in \NN$ and $c \geq 0$, 
the eigenvalues of \(\tilde{H}_c(m)\) come in plus/minus pairs
and $ (-2|c - 1|,2|c - 1|) \cap \sigma(\tilde{H}_c(m)) =\emptyset$.
This interval is the largest with this property.
\end{theorem}
In particular, $(-2|c - 1|,2|c - 1|)$ contains no spurious eigenvalues whatsoever.\\

Finally, we consider the `infinite dimensional limits',
that is, the matrices \(\mathbf{H}_c, \mathbf{A},
\mathbf{B}, \mathbf{T}_c,\mathbf{W}_c\), obtained from \(H, A,
B, T_c,W_c\) by stretching to infinity in both directions.
Thus $\bA=\bA^*$, $\bB=-\bB^*$, $\bT_c$, and 
$\bW_c=\bT_{-c}^*\bT_{-c}=\bW_c^*$ are bounded operators
on the Hilbert space $\ell^2(\ZZ)$, while $\bH_c$
and $\bK_c$ are selfadjoint operators on $\ell^2(\ZZ)\bigoplus \ell^2(\ZZ)$.
The operator 
\begin{equation*}
\bU =
\frac{1}{\sqrt{2}}
\left[\begin{array}{rr}
I     &       I \\
I     &      -I \\
\end{array}\right]
\end{equation*}
is unitary on  $\ell^2(\ZZ)\bigoplus \ell^2(\ZZ)$, 
where now $I$ denotes the identity on  $\ell^2(\ZZ)$.
The operators keep their algebraic relations
\begin{align*}
\bK_c = \bU^{-1}\bH_c\bU = \bU\bH_c\bU =
\left[\begin{array}{rr}
0    &  \bA_c - \bB \\
\bA_c + \bB   &  0 \\
\end{array}\right]
=
2
\left[\begin{array}{rr}
0    &  \bT_c \\
\bT_c^*  &  0 \\
\end{array}\right]. 
\end{align*}
Formula $\bW_c = \bT_{-c}^*\bT_{-c}$ gives us 
\[
(\mathbf{W}_c x)_j = -cx_{j-1} + (c^2 + 1)x_j - cx_{j+1}.
\]
Now using the isometric isomorphism \(\ell^2(\mathbb{Z}) \to L^2(0,2\pi)\) given by
\[
\psi(\lambda) = \frac{1}{\sqrt{2}}\sum_{k = -\infty}^\infty e^{ij\lambda} x_j,
\]
the operator \(\mathbf{W}_c\) goes over into
\begin{align*}
(\mathbf{X}_c\psi)(\lambda) &= 
\frac{1}{\sqrt{2}}\sum_{j = -\infty}^{\infty}e^{ij\lambda}(-cx_{j-1} + (c^2 + 1)x_j - cx_{j+1})
\\
&=\frac{1}{\sqrt{2}}\sum_{j = -\infty}^{\infty}e^{ij\lambda}
(c^2 + 1 -2c\cos\lambda)x_j
= (c^2 + 1 -2c \cos\lambda)\psi(\lambda), 
 \end{align*}
which is a multiplication operator with the spectrum

\begin{equation*}
 [-(1 + c)^2, -(1 - c)^2] \cup [(1 - c)^2, (1 + c)^2],
\end{equation*}
thus creating the spectral gap
\((|1 - c|, |1 + c|)\) of \(\mathbf{H}_c\). Thus, the obvious
approximation \(H_c\) obtained by cutting a 
`window' out of   \(\mathbf{H}_c\) gives rise to spectral pollution
in the spectral gap of $\bH_c$. By adding convenient boundary conditions
\(H_c\) we obtain the modification \(\tilde H_c\). This new approximation to 
\(\bH_c\) 
\begin{itemize}
\item[(i)] has no spectral pollution and 
\item[(ii)] keeps the symmetry
of the spectrum with respect to zero.
\end{itemize}

\paragraph*{Some numerical experiments.} 
Here we would like to report some interesting
observations based on numerical experiments. 
They are motivated by physical models of disordered systems.
In this context
the matrix \(A\) in (\ref{stolz})
is replaced by \(V_\omega=A + \diag(\omega_1,\ldots,\omega_n)\),
where \(\omega_i\) are independent random variables. 
Here we will consider a uniform distribution on the interval \([a,b]\). 
Then, as expected, the multiple 
eigenvalues \(\pm 2\) of the matrix 
\(H_\omega=H
\left[\begin{array}{rr}
A_\omega     &       B \\
-B    &      -A_\omega \\
\end{array}\right]
\) smear 
into uniformly distributed 
intervals,
but {\em the double small eigenvalue
is only slightly perturbed in the sense that for 
\(n\) large these two eigenvalues tend to zero}.
We illustrate this by exhibiting those eigenvalues of the matrix \(H_\omega\)
which are close to zero by taking \(a = -3,\ b = 3\).

\begin{verbatim}

        m  = 20              m = 50               m  = 100

       -7.2091e-01          -3.5965e-01          -0.170365
       -5.4659e-01          -1.4649e-01          -0.098804
       -5.4215e-04          -4.5522e-10          -9.819153e-31
        5.4215e-04           4.5522e-10           9.819153e-31
        5.4659e-01           1.4649e-01           0.098804
        7.2091e-01           3.5965e-01           0.170365


\end{verbatim}                 
We emphasize that the exhibited digits of 9.819153e-31 are accurate.
The phenomenon of two very small eigenvalues
is independent of the choice of \(a,b\). Note that the spurious eigenvalues
are not only small but about exponentially small as in the 
case of the constant diagonal, i.e.~$\omega_1=\dots=\omega_m= 2c$, studied above.\\

\begin{figure}[htp]
\begin{center}
\includegraphics[width=8cm]{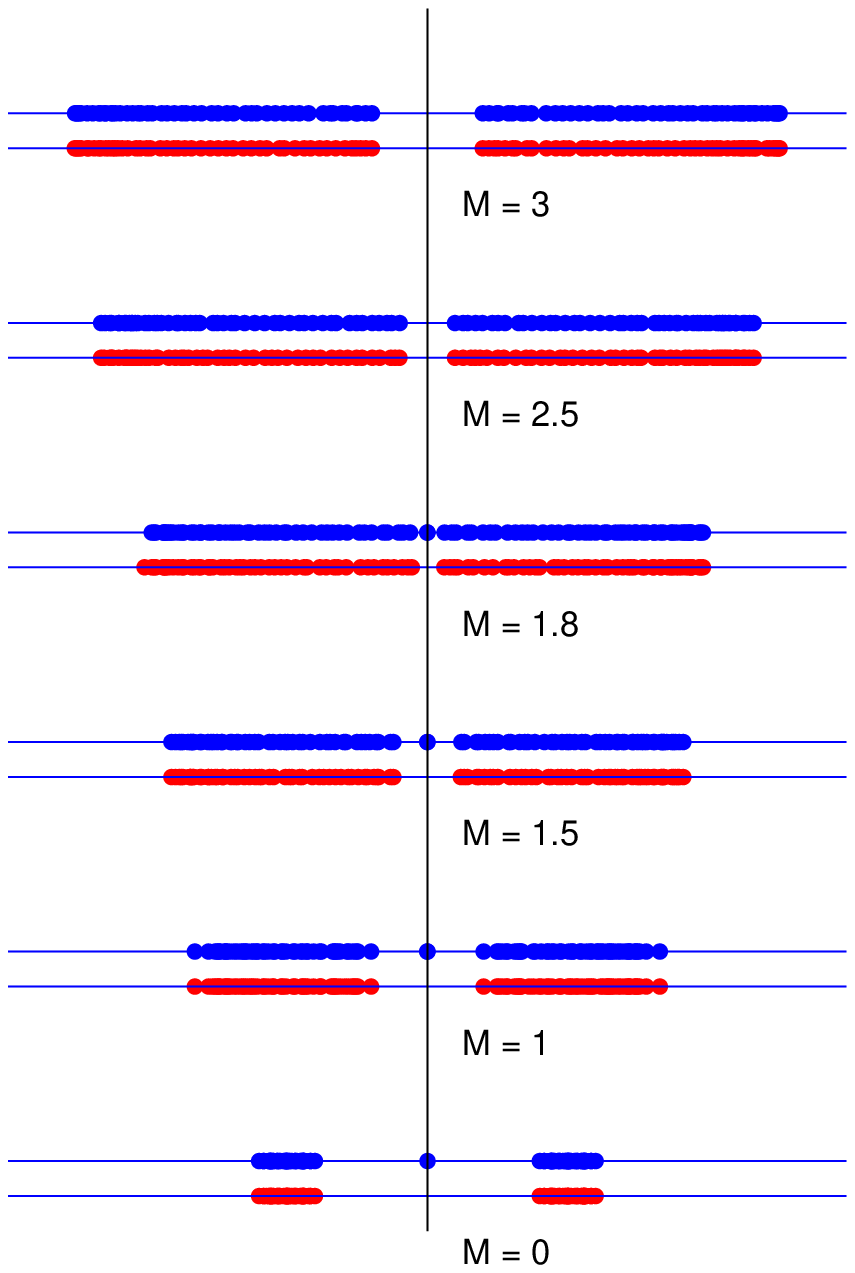}\caption{ \label{fig_0} Spectral gaps}
 \end{center}
 \end{figure}

Next we produce a series of graphically represented numerical results 
with
\[
[a,b] = [M - \delta ,M + \delta],
\]
\(n = 100\),  \(\delta = 0.5\), and \(M\) taking the values
\(0.1,\ 1,\ 1.5,\ 1.8,\ 2.5,\ 3\). They are contained in Figure \ref{fig_0}.

The upper line (blue) in a pair shows the spectrum of \(H_\omega\) 
and the lower (red) the one of
\(\tilde{H}_\omega\), obtained from \(H_\omega\) as in (\ref{HH}).\\ 

Summarising we may say:
{\em For \(M \approx 2\) or so the punctured spectral 
gap shrinks to zero;
then it starts growing again, but small eigenvalues are no more present,
because the matrix \(H_\omega\) has now become quasi-definite.\footnote{
Note that \(M\) as the mean value of \(\omega_i\) corresponds to the value
\(2c\) with \(H_c\) in (\ref{stolz}).} 
No theoretical explanation for the spurious small eigenvalues in this
case seems to
be available as yet. On the other hand, as expected, the matrix 
\(\tilde{H}_\omega\)
lacks the spurious small eigenvalues altogether. With \(\tilde{H}_\omega\)
there is no more symmetry of the spectrum with respect to zero.}\\


\begin{small}

{\bf Acknowledgement.} We are grateful to the colleagues who have 
obliged us with illuminating discussions and comments. These are
A.~B\"ottcher,
W.~Kirsch, 
T.~Lin{\ss},
I.~Naki\'{c}, 
M.~Omladi\v{c}, 
M.~Skrzipek,
and G.~Stolz.
\end{small}

\newpage
\appendix
\section{Some auxiliary computations}
\subsection*{Proof of equation \eqref{4x4pol}} 
Since the matrices of the type (\ref{4x4}) appear to be the source of
many illustrative examples we here give an explicit formula for
their eigenvalues (which come in plus/minus pairs).
We put
\[
A =
\left[\begin{array}{rr}
a_+     &    a \\
\bar{a}   &   a_-  \\
\end{array}\right],\quad
B =
\left[\begin{array}{rr}
b_+     &    b \\
\pm \bar{b}   &   b_-  \\
\end{array}\right],
\]
so that in the case of the minus sign in \(B\) the diagonal
elements \(b_\pm\) are purely imaginary. A straightforward
calculation gives
\[
\left(\lambda^2\right)_{1,2} =
\label{lambda4x4allg}
\frac{s}{2} \pm
\sqrt{s^2 - \left|a_+a_- - b_+b_- - |a|^2 + |b|^2\right|^2
 - \left|a_+b_- - b_+a_- - 2\Re\ \bar{a}b\right|^2}.
\]
with
\[
s = \frac{a_+^2 + a_-^2 + |b_+|^2 + |b_-|^2}{2} + |a|^2 + |b|^2.
\]
\subsection*{Proof of equation \eqref{kappa}} 
We derive the formula (\ref{kappa}).
Substituting \(x_j = A\cos j\alpha + B\sin j\alpha\) in (\ref{diffAhat})
we get
\begin{eqnarray*}
A\cos j\alpha\cos\alpha &+& A\sin j\alpha\sin\alpha
B\sin j\alpha\cos\alpha - B\cos j\alpha\sin\alpha\\ 
&+& \kappa(A\cos j\alpha + B\sin j\alpha) \\
+A\cos j\alpha\cos\alpha &-& A\sin j\alpha\sin\alpha
B\sin j\alpha\cos\alpha + B\cos j\alpha\sin\alpha\\
&=& 0
\end{eqnarray*}
or
\[
(A\cos j\alpha + B\sin j\alpha)(\cos\alpha + \cos\alpha + \kappa) = 0
\]
thus implying
\[
\kappa = -2\cos\alpha.
\]
The boundary conditions (\ref{bdAhat}) yield
\begin{eqnarray*}
A &=& -A\cos \alpha -B\sin \alpha,\\
A\cos (m+1)\alpha + B\sin (m+1)\alpha &=& A\cos m\alpha + B\sin m\alpha, 
\end{eqnarray*}
which is a homogeneous linear system
\begin{eqnarray*}
(1 + \cos\alpha)A + B\sin\alpha & = & 0,\\
A(\cos (m+1)\alpha - \cos m\alpha) + B(\sin (m+1)\alpha - \sin m\alpha)
 & = & 0, 
\end{eqnarray*}
so its determinant must vanish:
\[
(1 + \cos\alpha)2\cos\frac{(2m + 1)\alpha}{2}\sin\frac{\alpha}{2} +
2\sin\alpha\sin\frac{(2m + 1)\alpha}{2}\sin\frac{\alpha}{2} = 0,
\]
or equivalently,
\[
0 = \cos\frac{\alpha}{2}\cos\frac{(2m + 1)\alpha}{2} + 
\sin\frac{\alpha}{2}\sin\frac{(2m + 1)\alpha}{2} =
\cos\left(\frac{\alpha}{2} - \frac{(2m + 1)\alpha}{2}\right)
= \cos m\alpha.
\]
Hence
\[
\alpha = \alpha_k = \frac{2k + 1}{2m}\pi
\]
and
(\ref{kappa}) follows.
\end{document}